\pdfoutput=1
\RequirePackage{ifpdf}
\ifpdf 
\documentclass[pdftex]{sigma}
\else
\documentclass{sigma}
\fi

\usepackage[mathscr]{eucal}
\usepackage{enumerate}
\usepackage{xspace}
\usepackage[thinlines]{easymat}

\numberwithin{equation}{section}

\newtheorem{Theorem}{Theorem}[section]
\newtheorem{Lemma}[Theorem]{Lemma}

{\theoremstyle{definition}

\newtheorem{Remark}[Theorem]{Remark}
}

\newcommand{\vfld}[1]{\frac{\partial}{\partial{#1}}}

\DeclareMathOperator{\spn}{span}

\newcommand{\linsp}[1]{\spn{\left({#1}\right)}}

\begin{document}

\allowdisplaybreaks

\renewcommand{\PaperNumber}{034}

\FirstPageHeading

\ShortArticleName{Geometry of Optimal Control for Control-Af\/f\/ine Systems}

\ArticleName{Geometry of Optimal Control\\
for Control-Af\/f\/ine Systems}

\Author{Jeanne N.~CLELLAND~$^\dag$, Christopher G.~MOSELEY~$^\ddag$ and George R.~WILKENS~$^\S$}

\AuthorNameForHeading{J.N.~Clelland, C.G.~Moseley and G.R.~Wilkens}

\Address{$^\dag$~Department of Mathematics, 395 UCB, University of Colorado, Boulder, CO 80309-0395, USA}
\EmailD{\href{mailto:Jeanne.Clelland@colorado.edu}{Jeanne.Clelland@colorado.edu}}

\Address{$^\ddag$~Department of Mathematics and Statistics, Calvin College, Grand Rapids, MI 49546, USA}
\EmailD{\href{mailto:cgm3@calvin.edu}{cgm3@calvin.edu}}

\Address{$^\S$~Department of Mathematics, University of Hawaii at Manoa,\\
\hphantom{$^\S$}~2565 McCarthy Mall, Honolulu, HI 96822-2273, USA}
\EmailD{\href{mailto:grw@math.hawaii.edu}{grw@math.hawaii.edu}}

\ArticleDates{Received June 07, 2012, in f\/inal form April 03, 2013; Published online April 17, 2013}

\Abstract{Motivated by the ubiquity of control-af\/f\/ine systems in optimal control theory,
we investigate the geometry of
point-af\/f\/ine control systems with metric structures in dimensions two and three.
We compute local isometric invariants for point-af\/f\/ine distributions of constant type with metric structures for
systems with 2 states and 1 control and systems with 3 states and 1 control, and use Pontryagin's maximum principle to
f\/ind geodesic trajectories for homogeneous examples.
Even in these low dimensions, the behavior of these systems is surprisingly rich and varied.}


\Classification{58A30; 53C17; 58A15; 53C10}

\section{Introduction}

In~\cite{CMW09}, we investigated the local structure of \textit{point-affine distributions}.
A rank-\(s\) point-af\/f\/ine distribution on an $n$-dimensional manifold \(M\) is a~sub-bundle \({\EuScript F}\) of
the tangent bundle \(TM\) such that, for each \(x\in M\), the f\/iber \({\EuScript F}_{x} = T_{x} M \cap {\EuScript
F}\) is an \(s\)-dimensional af\/f\/ine subspace of~\(T_{x} M\) that contains a~distinguished point.
In local coordinates, the points of \({\EuScript F}\) are parametrized by \(s+1\) pointwise independent smooth vector
f\/ields \(v_{0}(x), v_{1}(x), \dotsc, v_{s}(x)\) for which \({\EuScript F}_{x} = v_{0}(x) + \linsp{v_{1}(x), \dotsc,
v_{s}(x)}\) and \(v_{0}(x)\) is the distinguished point in \({\EuScript F}_{x}\).

Our interest in point-af\/f\/ine distributions is motivated by a~family of ordinary dif\/ferential equations that
occurs in control theory: the control-af\/f\/ine systems.
A control system is a~system of underdetermined ODEs
\begin{gather*}
\dot{x}=f(x,u),
\end{gather*}
where $x\in M$ and $u$ takes values in an $s$-dimensional manifold ${\EuScript U}$.
The system is {\em control-affine} if the right-hand side is af\/f\/ine linear in the control variables $u$, i.e., if
the system locally has the form
\begin{gather}
\label{eq:control-system-with-drift}
\dot{x}(t)=v_{0}(x)+\sum_{i=1}^{s}v_{i}(x)u^{i}(t),
\end{gather}
where the controls \(u^{1}, \dotsc, u^{s}\) appear linearly in the right hand side and \(v_{0}, \dotsc, v_{s}\) are
\(s+1\) independent vector f\/ields (see, e.g.,~\cite{Jurdjevic97}). Replacing \(v_{0}\), which is called the {\em
drift} vector f\/ield, with a~linear combination of \(v_{1}, \dotsc, v_{s}\) added to \(v_{0}\) would yield an
equivalent system of dif\/ferential equations.
In many instances, however, there is a~distinguished null value for the controls (for example, consider turning of\/f
all motors on a~boat drifting downstream), and this null value determines a~distinguished drift vector f\/ield.
In these instances, we always choose \(v_{0}\) to be the distinguished drift vector f\/ield.
Consequently, the null value for the controls will be
\begin{gather*}
u^{1}=\dots=u^{s}=0.
\end{gather*}

While the control-af\/f\/ine systems~\eqref{eq:control-system-with-drift} may appear to be rather special, these
systems are ubiquitous.
In fact, any control system whatsoever becomes control-af\/f\/ine after a~single prolongation, so these systems
actually encompass all control systems, at the cost of increasing the number of state variables.

In~\cite{CMW09} we studied local dif\/feomorphism invariants for these point-af\/f\/ine structures.
A local equivalence for two point-af\/f\/ine structures is a~local dif\/feomorphism of $M$ whose derivative maps one
distinguished drift vector f\/ield to the other, and maps one af\/f\/ine sub-bundle to the other (see~\cite{CMW09} for
precise def\/initions).
With this notion of local equivalence, we were able to determine local normal forms for strictly af\/f\/ine, rank-\(1\)
point-af\/f\/ine structures of constant type when the manifold \(M\) had dimension 2 or 3.
In some cases the normal forms are parametrized by arbitrary functions.

The current paper seeks to ref\/ine the previous results by adding a~metric structure to the point-af\/f\/ine structure.
We do so by introducing a~positive def\/inite quadratic cost functional \(Q: {\EuScript F} \to \mathbb{R}\).
In local coordinates, where
\begin{displaymath}
w = v_{0}(x) + \sum_{i=1}^{s} v_{i}(x) u^{i} \in {\EuScript F}_{x},
\end{displaymath}
we will def\/ine
\begin{gather*}
Q_{x}(w)=\sum g_{i j}(x)u^{i}u^{j},
\end{gather*}
where the matrix \((g_{i j}(x))\) is positive def\/inite and the components are smooth functions of \(x\).
This is a~natural extension of the well-studied notion of a~sub-Riemannian metric on a~linear distribution, which
represents a~quadratic cost functional for a~driftless system (see, e.g.,~\cite{LS95,Montgomery02, Moseley01}).

With the added metric structure, we ref\/ine our notion of local point-af\/f\/ine equivalence to that of a~local
\textit{point-affine isometry}.
A local point-af\/f\/ine isometry is a~local point-af\/f\/ine equivalence that additionally preserves the quadratic
cost functional.

Let \( \gamma(t) = x(t)\) be a~trajectory for~\eqref{eq:control-system-with-drift}.
The added metric structure allows us to assign the following energy cost functional to \( \gamma(t)\):
\begin{gather}
\label{cost-functional-eq}
E(\gamma)=\frac{1}{2}\int_{\gamma}Q_{x(t)}\bigl(\dot{x}(t)\bigr)dt.
\end{gather}
Naturally associated to~\eqref{cost-functional-eq} is the \textit{optimal control problem} of f\/inding trajectories
of~\eqref{eq:control-system-with-drift} that minimize~\eqref{cost-functional-eq}.
We will use Pontryagin's maximum principle to f\/ind an ODE system on \(T^{*}M\) with the property that any minimal
cost trajectory for~\eqref{eq:control-system-with-drift} must be the projection of some solution for the ODE system on
\(T^{*}M\).

In this paper we shall only consider homogeneous examples, i.e., examples that admit a~symmetry group which acts
transitively on $M$.
We shall use the normal forms from~\cite{CMW09} as starting points, adding a~homogeneous metric structure to the
point-af\/f\/ine structure in each case.
Even in these low-dimensional cases, the analysis can be quite involved; we will see that these structures exhibit
surprisingly rich and varied behavior.

\section{Normal forms for homogeneous cases}

We begin by identifying the homogeneous examples of the point-af\/f\/ine systems described in~\cite{CMW09} where
possible, and then we describe the homogeneous metric structures on these systems.
In some cases, the metric structure must be added before the homogeneous examples can be identif\/ied.
Recall that the assumption of homogeneity is equivalent to the condition that all structure functions $T^i_{jk}$
appearing in the structure equations for a~canonical coframing are constants
(see~\cite{Gardner89} for details).

\subsection{Two states and one control}

In~\cite{CMW09}, we found two local normal forms under point-af\/f\/ine
equivalence.

{\bf Case 1.1.} \({\EuScript F} = \vfld{x^{1}} + \linsp{\vfld{x^{2}}}\).
The framing
\begin{gather*}
v_{1}=\vfld{x^{1}},
\qquad
v_{2}=\frac{1}{\lambda}\vfld{x^{2}}
\end{gather*}
(well-def\/ined up to scaling in \(v_{2}\)) has dual coframing
\begin{gather}
\eta^{1}=dx^{1},
\qquad
\eta^{2}=\lambda dx^{2},
\label{case-1-1-coframing}
\end{gather}
with structure equations
\begin{gather*}
d\eta^{1}=0,
\qquad
d\eta^{2}\equiv0
\mod\eta^{2}.
\end{gather*}
Because the method of equivalence does not lead to a~completely determined canonical cofra\-ming, it is not clear from
these structure equations whether this example is homogeneous as a~point-af\/f\/ine distribution.

Fortunately, this ambiguity is resolved when we add a~metric function to the point-af\/f\/ine structure.
This amounts to a~choice of function $G(x)>0$ for which the quadratic cost functional is given by
\begin{gather}
Q\left(\vfld{x^{1}}+u \vfld{x^{2}}\right)=\frac{1}{2}G(x)u^2.
\label{case-1-1-metric}
\end{gather}
For the point-af\/f\/ine structure, the frame vector~$v_2$ is only well-def\/ined up to a~scale factor; however, when
we impose a~metric structure~\eqref{case-1-1-metric}, we can choose~$v_2$ canonically (up to sign) by requiring that it
be a~unit vector for the metric.
This choice leads to a~{\em canonical} framing
\begin{gather*}
v_1=\vfld{x^1},
\qquad
v_2=\frac{1}{\sqrt{G(x)}}\vfld{x^2},
\end{gather*}
with corresponding canonical coframing
\begin{gather*}
\eta^1=dx^1,
\qquad
\eta^2=\sqrt{G(x)}\,dx^2.
\end{gather*}
The structure equations for this ref\/ined coframing are
\begin{gather*}
d\eta^{1}=0,
\qquad
d\eta^{2}=\frac{G_{x^{1}}}{2G}\eta^{1}\wedge\eta^{2},
\end{gather*}
and so the structure is homogeneous if and only if \(\frac{G_{x^{1}}}{2G}\) is equal to a~constant $c_1$.
This condition implies that{\samepage
\begin{gather*}
G\big(x^1,x^2\big)=G_{0}\big(x^{2}\big)e^{2c_1x^{1}}
\end{gather*}
for some function $G_0\big(x^2\big)$.}

The local coordinates in the coframing~\eqref{case-1-1-coframing} are only determined up to transformations of the form
\begin{gather}
x^1=\tilde{x}^1+a,
\qquad
x^2=\phi\big(\tilde{x}^2\big),
\label{case-1-1-coord-transformation}
\end{gather}
and under this transformation we have
\begin{gather*}
\tilde{G}_0\big(\tilde{x}^2\big)
=e^{2c_1a}\left(\big(\phi'\big(\tilde{x}^2\big)\right)^2G_0\big(\phi\big(\tilde{x}^2\big)\big).  
\end{gather*}
Therefore, we can apply a~transformation of the form~\eqref{case-1-1-coord-transformation} to arrange that
$\tilde{G}_0\big(\tilde{x}^2\big)=1$, and hence $\tilde{G}=e^{2c_1\tilde{x}^1}$.
Moreover, coordinates for which $G$ has this form are uniquely determined up to a~transformation of the form
\begin{gather*}
x^1=\tilde{x}^1+a,
\qquad
x^2=e^{-c_1a}\tilde{x}^2+b.
\end{gather*}

To summarize: the homogeneous metrics in this case are given by quadratic functionals of the form
\begin{gather*}
Q\left(\vfld{x^1}+u\vfld{x^2}\right)=\frac{1}{2}e^{2c_1x^1}u^2
\end{gather*}
for some constant $c_1$, with corresponding canonical coframings
\begin{gather*}
\eta^1=dx^1,
\qquad
\eta^2=e^{c_1x^1}\,dx^2.
\end{gather*}

{\bf Case 1.2.} \({\EuScript F} = x^{2} \Bigl(\vfld{x^{1}} + J\vfld{x^{2}}\Bigr) + \linsp{\vfld{x^{2}}}\).
We found a~canonical framing
\begin{gather}
v_{1}=x^{2}\Biggl(\vfld{x^{1}}+J\vfld{x^{2}}\Biggr),
\qquad
v_{2}=x^{2}\vfld{x^{2}},
\label{case-1-2-framing}
\end{gather}
with dual coframing
\begin{gather}
\eta^1=\frac{1}{x^{2}}dx^{1},
\qquad
\eta^{2}=\frac{1}{x^{2}}\big(dx^{2}-J dx^{1}\big),
\label{case-1-2-coframing}
\end{gather}
and structure equations
\begin{gather*}
d\eta^{1}=\eta^{1}\wedge\eta^{2},
\qquad
d\eta^{2}=T^{2}_{12}\eta^{1}\wedge\eta^{2},
\end{gather*}
where
\begin{gather}
T^{2}_{12}=x^{2}\frac{\partial J}{\partial x^{2}}-J.
\label{case-1-2-J-eqn-1}
\end{gather}
The structure is homogeneous if and only if $T^1_{12}$ is equal to a~constant $-j_0$.
According to equation~\eqref{case-1-2-J-eqn-1}, this is the case if and only if
\begin{gather}
J=x^2J_1\big(x^1\big)+j_0
\label{case-1-2-J-eqn}
\end{gather}
for some function $J_1\big(x^1\big)$.

The local coordinates in the coframing~\eqref{case-1-2-coframing} are only determined up to transformations of the form
\begin{gather}
x^{1}=\phi\big(\tilde{x}^{1}\big),
\qquad
x^{2}=\tilde{x}^{2}\phi'\big(\tilde{x}^{1}\big),
\label{case-1-2-coord-transformation}
\end{gather}
and under this transformation we have
\begin{gather*}
\tilde{J}\big(\tilde{x}^{1},\tilde{x}^{2}\big)=
J\bigl(\phi\big(\tilde{x}^{1}\big),\tilde{x}^{2}\phi'\big(\tilde{x}^{1}\big)\bigr)
-\tilde{x}^{2}\frac{\phi''\big(\tilde{x}^{1}\big)}{\phi'\big(\tilde{x}^{1}\big)}.
\end{gather*}
In the homogeneous case~\eqref{case-1-2-J-eqn}, this implies that
\begin{gather*}
\tilde{J}_1\big(\tilde{x}^{1}\big)
=\phi'\big(\tilde{x}^{1}\big)J_1\big(\phi\big(\tilde{x}^{1}\big)\big)
-\frac{\phi''\big(\tilde{x}^{1}\big)}{\phi'\big(\tilde{x}^{1}\big)}.
\end{gather*}
Therefore, we can apply a~transformation of the form~\eqref{case-1-2-coord-transformation} to arrange that
$\tilde{J}_1\big(\tilde{x}^{1}\big)=0$, and hence $\tilde{J}=j_0$.
Moreover, coordinates for which $J$ is constant are uniquely determined up to an af\/f\/ine transformation
\[
x^{1} = a \tilde{x}^{1} + b,
\qquad
x^{2} = a \tilde{x}^{2}.
\]

Now suppose that a~metric on the point-af\/f\/ine structure is given by
\begin{gather}
Q\left(v_1+u v_2\right)=Q\left(x^2\left(\vfld{x^{1}}+j_0\vfld{x^2}\right)+u\left(x^2 \vfld{x^{2}}\right) \right)=
\frac{1}{2}G(x)u^2.
\label{case-1-2-metric}
\end{gather}
This case dif\/fers from the previous case in that the control vector f\/ield $v_2$ is already canonically def\/ined by
the point-af\/f\/ine structure prior to the introduction of a~metric.
Therefore, in order that the metric~\eqref{case-1-2-metric} be homogeneous, the unit control vector f\/ield
\begin{gather*}
\frac{1}{\sqrt{G(x)}} v_2
\end{gather*}
must be a~{\em constant} scalar multiple of $v_2$.
Thus we must have $G(x) = g_0$ for some positive constant $g_0$, and the homogeneous metrics in this case are given by
quadratic functionals of the form
\begin{gather*}
Q(v_1+u v_2)=\frac{1}{2}g_0u^2
\end{gather*}
for some positive constant $g_0$, where $v_1$, $v_2$ are the canonical frame vectors~\eqref{case-1-2-framing}.

\subsection{Three states and one control}

In~\cite{CMW09}, we found three nontrivial local normal forms under point-af\/f\/ine equivalence.
\begin{Remark}
This classif\/ication assumes that the point-af\/f\/ine distribution is either bracket-generating or almost
bracket-generating; otherwise the 3-manifold $M$ can locally be foliated by a~1-parameter family of 2-dimensional
submanifolds such that every trajectory of ${\EuScript F}$ is contained in a~single leaf of the foliation.
\end{Remark}

{\bf Case 2.1.} \({\EuScript F} = \left(\vfld{x^{1}} + x^3 \vfld{x^2} + J \vfld{x^3} \right) +
\linsp{\vfld{x^{3}}}\).
The framing
\begin{gather*}
v_{1}=\vfld{x^{1}}+x^3\vfld{x^2}+J\vfld{x^3},
\qquad
v_{2}=\vfld{x^{3}},
\qquad
v_3=-[v_1,v_2]=\vfld{x^2}+J_{x^3}\vfld{x^3}
\end{gather*}
(well-def\/ined up to dilation in the \((v_{2}, v_3)\)-plane) has dual coframing
\begin{gather*}
\eta^{1}=dx^{1},
\qquad
\eta^{2}=dx^3-J\,dx^1-J_{x^3}\big(dx^2-x^3\,dx^1\big),
\qquad
\eta^3=dx^2-x^3\,dx^1,
\end{gather*}
with structure equations
\begin{gather*}
\begin{split}
& d\eta^{1} = 0,
\\
& d\eta^{2}  \equiv T^2_{13} \eta^1 \wedge \eta^3
\mod \eta^2,
\\
& d\eta^3 \equiv \eta^1 \wedge \eta^2
\mod \eta^3.
\end{split}
\end{gather*}
As in Case~1.1, the method of equivalence does not lead to a~completely determined coframing, so it is not clear from
these structure equations whether this example is homogeneous as a~point-af\/f\/ine distribution.

So, suppose that a~metric on the point-af\/f\/ine structure is given by
\begin{gather}
Q\left(\left(\vfld{x^{1}}+x^3\vfld{x^2}+J\vfld{x^3}\right)+u\,\vfld{x^{3}}\right)=\frac{1}{2}G(x)u^2.
\label{case-2-1-metric}
\end{gather}
The addition of the metric~\eqref{case-2-1-metric} allows us to choose a~canonical framing (up to sign) by requi\-ring~$v_2$ to be a~unit vector for the metric, i.e.,
\begin{gather*}
v_2=\frac{1}{\sqrt{G(x)}} \vfld{x^3},
\end{gather*}
and setting
\begin{gather*}
v_3=-[v_1,v_2].
\end{gather*}
The canonical coframing associated to this framing is given by
\begin{gather}
\eta^{1}=dx^{1},
\qquad\!
\eta^{2}\equiv\sqrt{G(x)}\big(dx^3-J\,dx^1\big)\mod{\eta^3},
\qquad\!
\eta^3=\sqrt{G(x)}\big(dx^2-x^3\,dx^1\big).\!\!\!
\label{case-2-1-local-coord-exp}
\end{gather}

In order to identify the homogeneous examples, we consider the structure equations for the
coframing~\eqref{case-2-1-local-coord-exp}, taking into account the fact that local coordinates for which the coframing
takes the form~\eqref{case-2-1-local-coord-exp} are determined only up to transformations of the form
\begin{gather}
\label{case-2-1-local-coord-trans}
x^1=\tilde{x}^1+a,
\qquad
x^2=\phi\big(\tilde{x}^1,\tilde{x}^2\big),
\qquad
x^3=\phi_{\tilde{x}^1}\big(\tilde{x}^1,\tilde{x}^2\big)+\tilde{x}^3\phi_{\tilde{x}^2}\big(\tilde{x}^1,\tilde{x}^2\big),
\end{gather}
with $\phi_{\tilde{x}^2}\neq0$.
Under such a~transformation we have
\begin{gather}
\sqrt{\tilde{G}\big(\tilde{x}^1,\tilde{x}^2,\tilde{x}^3\big)}=\sqrt{G\big(x^1,x^2,x^3\big)} \phi_{\tilde{x}^2},
\label{case-2-1-new-G}
\\
\tilde{J}\big(\tilde{x}^1,\tilde{x}^2,\tilde{x}^3\big)=
\frac{1}{\phi_{\tilde{x}^2}}\left(J\big(x^1,x^2,x^3\big)-\phi_{\tilde{x}^2\tilde{x}^2}
\big(\tilde{x}^3\big)^2-2\phi_{\tilde{x}^1\tilde{x}^2}\tilde{x}^3-\phi_{\tilde{x}^1\tilde{x}^1}\right),
\label{case-2-1-new-J}
\end{gather}
with $x^1$, $x^2$, $x^3$ as in~\eqref{case-2-1-local-coord-trans}.

First consider the structure equation for $d\eta^3$.
A computation shows that
\begin{gather*}
d\eta^3\equiv\frac{G_{x^3}}{2G^{{3/2}}}\eta^2\wedge\eta^3\mod{\eta^1}.
\end{gather*}
Therefore, homogeneity implies that $\frac{G_{x^3}}{2G^{{3/2}}}$ must be equal to a~constant $-c_1$.
The remaining analysis varies considerably depending on whether $c_1$ is zero or nonzero.

{\bf Case 2.1.1.}
First suppose that $c_1=0$.
Then $G_{x^3}=0$, and so
\begin{gather*}
G\big(x^1,x^2,x^3\big)=G_0\big(x^1,x^2\big)
\end{gather*}
for some function $G_0\big(x^1,x^2\big)$.
According to~\eqref{case-2-1-new-G}, by a~local change of coordinates of the form~\eqref{case-2-1-local-coord-trans}
with $\phi$ a~solution of the PDE
\begin{gather*}
\phi_{\tilde{x}^2}\big(\tilde{x}^1,\tilde{x}^2\big)
=\frac{1}{G_0\big(\tilde{x}^1,\phi\big(\tilde{x}^1,\tilde{x}^2\big)\big)},
\end{gather*}
we can arrange that $\tilde{G}_0\big(\tilde{x}^1,\tilde{x}^2\big)=1$.
This condition is preserved by transformations of the form~\eqref{case-2-1-local-coord-trans} with
\begin{gather}
\phi\big(\tilde{x}^1,\tilde{x}^2\big)=\tilde{x}^2+\phi_0\big(\tilde{x}^{1}\big).
\label{case-2-1-1-local-coord-trans-refined}
\end{gather}

With the assumption that $G\big(x^1,x^2,x^3\big)=1$, the equation for $d\eta^3$ reduces to
\begin{gather*}
d\eta^3=\eta^1\wedge\eta^2+J_{x^3}\eta^1\wedge\eta^3.
\end{gather*}
Therefore, $J_{x^3}$ must be equal to a~constant $c_3$, and so
\begin{gather*}
J\big(x^1,x^2,x^3\big)=c_3x^3+J_0\big(x^1,x^2\big)
\end{gather*}
for some function $J_0\big(x^1,x^2\big)$.
Now the equation for $d\eta^2$ becomes
\begin{gather*}
d\eta^2=(J_0)_{x^2}\,\eta^1\wedge\eta^3.
\end{gather*}
Therefore, $(J_0)_{x^2}$ must be equal to a~constant $c_2$, and so
\begin{gather*}
J_0\big(x^1,x^2\big)=c_2x^2+J_1\big(x^1\big)
\end{gather*}
for some function $J_1\big(x^1\big)$.
With $\phi$ as in~\eqref{case-2-1-1-local-coord-trans-refined} and
\begin{gather*}
J\big(x^1,x^2,x^3\big)=c_2x^2+c_3x^3+J_1\big(x^1\big),
\end{gather*}
equation~\eqref{case-2-1-new-J} reduces to
\begin{gather*}
\tilde{J}_1\big(\tilde{x}^1\big)=
J_1\big(\tilde{x}^1+a\big)-\left(\phi_0''\big(\tilde{x}^{1}\big)-c_3\phi_0'\big(\tilde{x}^{1}\big)
-c_2\phi_0\big(\tilde{x}^{1}\big)\right).
\end{gather*}
Therefore, we can choose local coordinates to arrange that $\tilde{J}_1\big(\tilde{x}^{1}\big)=0$.

To summarize, we have constructed local coordinates for which
\begin{gather*}
G\big(x^1,x^2,x^3\big)=1,
\qquad
J\big(x^1,x^2,x^3\big)=c_2x^2+c_3x^3.
\end{gather*}
These coordinates are determined up to transformations of the form
\begin{gather*}
x^1=\tilde{x}^1+a,
\qquad
x^2=\tilde{x}^2+\phi_0\big(\tilde{x}^{1}\big),
\qquad
x^3=\tilde{x}^3+\phi_0'\big(\tilde{x}^{1}\big),
\end{gather*}
where $\phi_0\big(\tilde{x}^{1}\big)$ is a~solution of the ODE
\begin{gather*}
\phi_0''\big(\tilde{x}^1\big)-c_3\phi_0'\big(\tilde{x}^{1}\big)-c_2\phi_0\big(\tilde{x}^{1}\big)=0.
\end{gather*}

{\bf Case 2.1.2.} Now suppose that $c_1\neq0$.
Then
\begin{gather*}
G\big(x^1,x^2,x^3\big)=\frac{1}{\left(c_1x^3+G_0\big(x^1,x^2\big)\right)^2}
\end{gather*}
for some function $G_0\big(x^1,x^2\big)$.
According to~\eqref{case-2-1-new-G}, by a~local change of coordinates of the form~\eqref{case-2-1-local-coord-trans}
with $\phi$ a~solution of the PDE
\begin{gather*}
\phi_{x^1}\big(\tilde{x}^1,\tilde{x}^2\big)
=\frac{1}{c_1}G_0\left(\tilde{x}^1,\phi\big(\tilde{x}^1,\tilde{x}^2\big)\right),
\end{gather*}
we can arrange that $\tilde{G}_0\big(\tilde{x}^1,\tilde{x}^2\big)=0$.
This condition is preserved by transformations of the form~\eqref{case-2-1-local-coord-trans} with
\begin{gather}
\phi\big(\tilde{x}^1,\tilde{x}^2\big)=\phi_0\big(\tilde{x}^2\big).
\label{case-2-1-2-local-coord-trans-refined}
\end{gather}

With the assumption that $G\big(x^1,x^2,x^3\big)=\frac{1}{(c_1x^3)^2}$, the equation for $d\eta^3$ reduces to
\begin{gather*}
d\eta^3=\eta^1\wedge\eta^2-\frac{(2J-x^3J_{x^3})}{x^3}\eta^1\wedge\eta^3-c_1\eta^2\wedge\eta^3.
\end{gather*}
Therefore, $\frac{(2J-x^3J_{x^3})}{x^3}$ must be equal to a~constant $c_3$, and so
\begin{gather*}
J\big(x^1,x^2,x^3\big)=c_3x^3+J_0\big(x^1,x^2\big)\big(x^3\big)^2
\end{gather*}
for some function $J_0\big(x^1,x^2\big)$.
Now the equation for $d\eta^2$ becomes
\begin{gather*}
d\eta^2=-x^3(J_0)_{x^1}\,\eta^1\wedge\eta^3.
\end{gather*}
The quantity $-x^3(J_0)_{x^1}$ can only be constant if $(J_0)_{x^1}=0$; therefore, we must have
\begin{gather*}
J_0\big(x^1,x^2\big)=J_1\big(x^2\big)
\end{gather*}
for some function $J_1\big(x^2\big)$.
With $\phi$ as in~\eqref{case-2-1-2-local-coord-trans-refined} and
\begin{gather*}
J\big(x^1,x^2,x^3\big)=c_3x^3+J_1\big(x^2\big)\big(x^3\big)^2,
\end{gather*}
equation~\eqref{case-2-1-new-J} reduces to
\begin{gather*}
\tilde{J}_1\big(\tilde{x}^2\big)
=J_1\big(\phi_0\big(\tilde{x}^2\big)\big)\phi_0'\big(\tilde{x}^2\big)
-\frac{\phi_0''\big(\tilde{x}^2\big)}{\phi_0'\big(\tilde{x}^2\big)}.
\end{gather*}
Therefore, we can choose local coordinates to arrange that $\tilde{J}_1\big(\tilde{x}^2\big)=0$.

To summarize, we have constructed local coordinates for which
\begin{gather*}
G\big(x^1,x^2,x^3\big)=\frac{1}{\big(c_1x^3\big)^2},
\qquad
J\big(x^1,x^2,x^3\big)=c_3x^3.
\end{gather*}
These coordinates are determined up to transformations of the form
\begin{gather*}
x^1=\tilde{x}^1+a,
\qquad
x^2=b\tilde{x}^2+c,
\qquad
x^3=b\tilde{x}^3+c.
\end{gather*}

{\bf Case 2.2.}
\({\EuScript F} = \left(x^2 \vfld{x^{1}} + x^3 \vfld{x^2} + J \left(x^2
\vfld{x^3}\right) \right) + \linsp{\vfld{x^{3}}}\).
We found a~canonical framing
\begin{gather}
v_{1}=x^{2}\vfld{x^{1}}+x^3\vfld{x^2}+J\left(x^2\vfld{x^{3}}\right),
\nonumber
\\
v_{2}=x^{2}\vfld{x^{3}},
\nonumber
\\
v_3=-[v_1,v_2]=x^2\vfld{x^2}+\left(\big(x^2\big)^2J_{x^3}-x^3\right)\vfld{x^3},
\label{case-2-2-framing}
\end{gather}
with dual coframing
\begin{gather}
\eta^1=\frac{1}{x^2}dx^1,
\nonumber
\\
\eta^2=
\frac{1}{x^2}dx^3-\frac{1}{x^2}J\,dx^1-\left(J_{x^3}-\frac{x^3}{\big(x^2\big)^2}\right)
\left(dx^2-\frac{x^3}{x^2}dx^1\right),
\nonumber
\\
\eta^3=\frac{1}{x^2}dx^2-\frac{x^3}{\big(x^2\big)^2}dx^1,
\label{case-2-2-coframing}
\end{gather}
and structure equations
\begin{gather}
d\eta^{1}=\eta^{1}\wedge\eta^{3},
\nonumber\\
d\eta^{2}=T^{2}_{13}\eta^{1}\wedge\eta^{3}+T^{2}_{23}\eta^{2}\wedge\eta^{3},
\nonumber\\
d\eta^3=\eta^1\wedge\eta^2+T^3_{13}\eta^1\wedge\eta^3.\label{case-2-2-structure-eqs}
\end{gather}

The local coordinates in the coframing~\eqref{case-2-2-coframing} are only determined up to transformations of the form
\begin{gather}
x^1=\phi\big(\tilde{x}^{1}\big),
\qquad
x^2=\phi'\big(\tilde{x}^{1}\big)\tilde{x}^2,
\qquad
x^3=\phi'\big(\tilde{x}^{1}\big)\tilde{x}^3+\phi''\big(\tilde{x}^{1}\big)\big(\tilde{x}^2\big)^2,
\label{case-2-2-local-coord-trans}
\end{gather}
with $\phi'\big(\tilde{x}^{1}\big)\neq0$.
Under such a~transformation we have
\begin{gather}
\tilde{J}\big(\tilde{x}^1,\tilde{x}^2,\tilde{x}^3\big)=
J\big(x^1,x^2,x^3\big)-\frac{1}{\phi'\big(\tilde{x}^{1}\big)}
\left(\phi'''\big(\tilde{x}^{1}\big)\big(\tilde{x}^2\big)^2+3\phi''\big(\tilde{x}^{1}\big)\tilde{x}^3\right),
\label{case-2-2-new-J}
\end{gather}
with $x^1$, $x^2$, $x^3$ as in~\eqref{case-2-2-local-coord-trans}.

First consider the structure equation for~$\eta^3$.
Substituting the expressions~\eqref{case-2-2-coframing} into the structure equation~\eqref{case-2-2-structure-eqs} for
$d\eta^3$ shows that
\begin{gather*}
T^2_{12}=x^2J_{x^3}-3\frac{x^3}{x^2}.
\end{gather*}
Homogeneity implies that $T^2_{12}$ must be equal to a~constant $a$, from which it follows that
\begin{gather*}
J\big(x^1,x^2,x^3\big)=\frac{3}{2}\left(\frac{x^3}{x^2}\right)^2+a\frac{x^3}{x^2}+J_0\big(x^1,x^2\big)
\end{gather*}
for some function $J_0\big(x^1,x^2\big)$.
Now the equation for $d\eta^2$ yields
\begin{gather*}
T^2_{13}=x^2(J_0)_{x^2}-2J_0-a\frac{x^3}{x^2},
\end{gather*}
and homogeneity implies that $T^2_{13}$ must be constant.
The quantity $\big(x^2(J_0)_{x^2}-2J_0-a \frac{x^3}{x^2}\big)$ can only be constant if $a=0$; therefore, we must
have $a=0$ and
\begin{gather*}
x^2(J_0)_{x^2}-2J_0=-2c_1
\end{gather*}
for some constant $c_1$.
Therefore,
\begin{gather*}
J_0\big(x^1,x^2\big)=c_1+J_1\big(x^1\big)\big(x^2\big)^2
\end{gather*}
for some function $J_1\big(x^1\big)$, and
\begin{gather*}
J\big(x^1,x^2,x^3\big)=\frac{3}{2}\left(\frac{x^3}{x^2}\right)^2+c_1+J_1\big(x^1\big)\big(x^2\big)^2.
\end{gather*}
With $\phi$ as in~\eqref{case-2-2-local-coord-trans} and $J$ as above, equation~\eqref{case-2-2-new-J} reduces to
\begin{gather*}
\tilde{J}_1\big(\tilde{x}^{1}\big)=
\phi'\big(\tilde{x}^{1}\big)^2J_1\big(\phi\big(\tilde{x}^{1}\big)\big)
-\frac{\phi'''\big(\tilde{x}^{1}\big)}{\phi'\big(\tilde{x}^{1}\big)}
+\frac{3}{2}\frac{\phi''\big(\tilde{x}^{1}\big)}{\big(\phi'\big(\tilde{x}^{1}\big)\big)^2}.
\end{gather*}
Therefore, we can choose local coordinates to arrange that $\tilde{J}_1\big(\tilde{x}^{1}\big)=0$.
This condition is pre\-served by transformations of the form~\eqref{case-2-2-local-coord-trans} with
\begin{gather*}
\frac{\phi'''\big(\tilde{x}^{1}\big)}{\phi'\big(\tilde{x}^{1}\big)}
-\frac{3}{2}\frac{\phi''\big(\tilde{x}^{1}\big)}{\big(\phi'\big(\tilde{x}^{1}\big)\big)^2}=0.
\end{gather*}
This implies that $\phi$ is a~linear fractional transformation, i.e.,
\begin{gather*}
\phi\big(\tilde{x}^{1}\big)=\frac{a\tilde{x}^1+b}{c\tilde{x}^1+d}.
\end{gather*}

Now suppose that a~metric on the point-af\/f\/ine structure is given by
\begin{gather}
Q\left(v_1+u v_2\right)=\frac{1}{2}G(x)u^2.
\label{case-2-2-metric}
\end{gather}
As in Case 1.2, the control vector f\/ield $v_2$ is already canonically def\/ined by the point-af\/f\/ine structure
prior to the introduction of a~metric.
Therefore, in order that the metric~\eqref{case-2-2-metric} be homogeneous, the unit control vector f\/ield
\begin{gather*}
\frac{1}{\sqrt{G(x)}}\,v_2
\end{gather*}
must \looseness=-1 be a~constant scalar multiple of $v_2$.
Thus we must have $G(x)=g_0$ for some positive cons\-tant~$g_0$, and the homogeneous metrics in this case are given by
quadratic functionals of the form
\begin{gather*}
Q(v_1+u v_2)=\frac{1}{2}g_0u^2
\end{gather*}
for some positive constant $g_0$, where $v_1$, $v_2$, $v_3$ are the canonical frame vectors~\eqref{case-2-2-framing}.

To summarize, we have constructed local coordinates for which
\begin{gather*}
G\big(x^1,x^2,x^3\big)=g_0,
\qquad
J\big(x^1,x^2,x^3\big)=\frac{3}{2}\left(\frac{x^3}{x^2}\right)^2+c_1.
\end{gather*}
These coordinates are determined up to transformations of the form
\begin{gather*}
x^1=\frac{a\tilde{x}^1+b}{c\tilde{x}^1+d},
\qquad
x^2=\frac{ad-bc}{(c\tilde{x}^1+d)^2}\tilde{x}^2,
\qquad
x^3=\frac{ad-bc}{(c\tilde{x}^1+d)^2}\tilde{x}^3-\frac{2c(ad-bc)}{(c\tilde{x}^1+d)^3}\tilde{x}^2.
\end{gather*}

{\bf Case 2.3.}
\begin{gather*}
{\EuScript F}=\left(\vfld{x^{1}}+J\left(x^3\vfld{x^1}+\vfld{x^2}+H\vfld{x^3}\right)\right)     
+\linsp{x^3\vfld{x^1}+\vfld{x^2}+H\vfld{x^3}},
\end{gather*}
where $\frac{\partial H}{\partial x^1}\neq0$.
We found a~canonical framing
\begin{gather*}
v_{1}=\vfld{x^{1}}+J\left(x^3\vfld{x^1}+\vfld{x^2}+H\vfld{x^3}\right),
\\
v_{2}=\frac{\epsilon}{\sqrt{\epsilon H_{x^1}}}\left(x^3\vfld{x^1}+\vfld{x^2}+H\vfld{x^3}\right),
\\
v_3=-[v_1,v_2],
\end{gather*}
where $\epsilon=\pm1=\text{sgn}(H_{x^1})$, with dual coframing{\samepage
\begin{gather}
\eta^{1}=dx^{1}-x^3\,dx^2,
\nonumber\\
\eta^2\equiv\epsilon\sqrt{\epsilon H_{x^1}}\left(dx^2-J\big(dx^1-x^3\,dx^2\big)\right)\mod{\eta^3},
\nonumber
\\
\eta^3=\frac{1}{\sqrt{\epsilon H_{x^1}}}\left(H\,dx^2-dx^3\right),
\label{case-2-3-coframing}
\end{gather}
and structure equations}
\begin{gather}
d\eta^1=T^1_{13}\eta^1\wedge\eta^3+T^1_{23}\eta^2\wedge\eta^3,
\nonumber\\
d\eta^2=T^2_{13}\eta^1\wedge\eta^3+T^2_{23}\eta^2\wedge\eta^3,
\nonumber
\\
d\eta^3=\eta^1\wedge\eta^2+T^3_{13}\eta^1\wedge\eta^3+T^3_{23}\eta^2\wedge\eta^3.
\label{case-2-3-structure-eqs}
\end{gather}
The identif\/ication of homogeneous examples is considerably more complicated than in the previous cases.
We refer the reader to Appendix~\ref{details-app} for the details.
We f\/ind that the homogeneous examples in this case are all locally equivalent to one of the following:
\begin{itemize}\itemsep=0pt
\item $J\big(x^1,x^2,x^3\big)=c_1$, $H\big(x^1,x^2,x^3\big)=\epsilon\big(x^1+c_2x^3\big)$
\\
 for some constants $c_1$,
$c_2$;
\item $J\big(x^1,x^2,x^3\big)=c_1\cos\big(c_3x^1\big)\big/\sqrt{\epsilon c_3\big(c_3\big(x^3\big)^2+c_4\big)}$,
\\
 $H\big(x^1,x^2,x^3\big)
=\big(c_3\big(x^3\big)^2+c_4\big)\tan\big(c_3x^1\big)+F_{20}\big(x^2\big)\sqrt{c_3\big(x^3\big)^2+c_4}$
\\
 for some constants $c_1$, $c_3$, $c_4$ with $c_3\neq0$, and some arbitrary function $F_{20}\big(x^2\big)$;
 \item
$J\big(x^1,x^2,x^3\big)=c_1\cosh\big(c_3x^1\big)\big/\sqrt{\epsilon c_3(c_3\big(x^3\big)^2-c_4)}$,
\\
 $H\big(x^1,x^2,x^3\big)
=\big({-}c_3\big(x^3\big)^2+c_4\big)\tanh\big(c_3x^1\big)+F_{20}\big(x^2\big)\sqrt{c_3\big(x^3\big)^2-c_4}$
\\
 for some constants $c_1$, $c_3$, $c_4$ with $c_3\neq0$, and some arbitrary function $F_{20}\big(x^2\big)$.
\end{itemize}

Now suppose that a~metric on the point-af\/f\/ine structure is given by
\begin{gather*}
Q\left(v_1+u v_2\right)=\frac{1}{2}G(x)u^2.
\end{gather*}
As in the previous case, since the control vector f\/ield $v_2$ is already canonically def\/ined by the
point-af\/f\/ine structure prior to the introduction of a~metric, we must have $G(x)=g_0$ for some positive constant
$g_0$.

The results of this section are encapsulated in the following two theorems:
\begin{Theorem}
\label{2-states-1-control-homog-examples-thm}
Let ${\EuScript F}$ be a~rank $1$ strictly affine point-affine distribution of constant type on a~$2$-dimensional
manifold $M$, equipped with a~positive definite quadratic cost functional~$Q$.
If the structure $({\EuScript F},Q)$ is homogeneous, then $({\EuScript F},Q)$ is locally point-affine equivalent to
\begin{gather*}
{\EuScript F}=v_1+\linsp{v_2},
\qquad
Q(v_1+u v_2)=\frac{1}{2}G(x)u^2,
\end{gather*}
where the triple $(v_1,v_2,G(x))$ is one of the following:
\begin{alignat*}{5}
&  (1.1) \quad &&  v_1=\vfld{x^1},
\qquad &&
v_2=\vfld{x^2},
\qquad &&
G(x)=e^{2c_1x^1}; &\\
& (1.2)\quad && v_1=x^{2}\left(\vfld{x^{1}}+j_0\vfld{x^{2}}\right),
\qquad &&
v_2=x^2\vfld{x^{2}},
\qquad &&
G(x)=g_0. &
\end{alignat*}
\end{Theorem}

\begin{Theorem}
\label{3-states-1-control-homog-examples-thm}
Let ${\EuScript F}$ be a~rank $1$, strictly affine, bracket-generating or almost bracket-generating point-affine
distribution of constant type on a~$3$-dimensional manifold $M$, equipped with a~positive definite quadratic cost
functional~$Q$.
If the structure $({\EuScript F},Q)$ is homogeneous, then $({\EuScript F},Q)$ is locally point-affine equivalent to
\begin{gather*}
{\EuScript F}=v_1+\linsp{v_2},
\qquad
Q(v_1+u v_2)=\frac{1}{2}G(x)u^2,
\end{gather*}
where the triple $(v_1,v_2,G(x))$ is one of the following:
\begin{alignat*}{3}
& (2.1.1) \quad && v_1=\vfld{x^{1}}+x^3\vfld{x^2}+\big(c_2x^2+c_3x^3\big)\vfld{x^3},
\qquad
v_2=\vfld{x^3},
\qquad
G(x)=1; & \\
& (2.1.2) \quad && v_1=\vfld{x^{1}}+x^3\vfld{x^2}+c_3x^3\vfld{x^3},
\qquad
v_2=\vfld{x^3},
\qquad
G(x)=\frac{1}{\big(c_1x^3\big)^2}; &
\\
&(2.2) \quad &&
v_1=x^{2}\vfld{x^{1}}+x^3\vfld{x^2}
+\left(\frac{3}{2}\left(\frac{x^3}{x^2}\right)^2+c_1\right)\left(x^2\vfld{x^{3}}\right), &\\
&&& v_2=x^{2}\vfld{x^{3}},
\qquad
G(x)=g_0; &
\\
& (2.3.1) \quad && v_1=\vfld{x^{1}}+c_1\left(x^3\vfld{x^1}+\vfld{x^2}+\epsilon\big(x^1+c_2x^3\big)\vfld{x^3}\right), &
\\
&&& v_2=\epsilon\left(x^3\vfld{x^1}+\vfld{x^2}+\epsilon\big(x^1+c_2x^3\big)\vfld{x^3}\right),
\qquad
G(x)=g_0; &
\\
& (2.3.2) \quad && v_1=\vfld{x^{1}}+\frac{c_1\cos(c_3x^1)}
{\sqrt{\epsilon c_3(c_3(x^3)^2+c_4)}}
\left(x^3\vfld{x^1}+\vfld{x^2}+H\vfld{x^3}\right), &
\\
&&& v_2=\epsilon\left(x^3\vfld{x^1}+\vfld{x^2}+H\vfld{x^3}\right),
\qquad
G(x)=g_0, &
\\
&&& \text{where} \quad H=\left(\big(c_3\big(x^3\big)^2+c_4\big)\tan\big(c_3x^1\big)+F_{20}\big(x^2\big)\sqrt{c_3\big(x^3\big)^2+c_4}\right); &
\\
& (2.3.3) \quad && v_1=\vfld{x^{1}}+\frac{c_1\cosh(c_3x^1)}
{\sqrt{\epsilon c_3(c_3(x^3)^2-c_4)}}\left(x^3\vfld{x^1}+\vfld{x^2}+H\vfld{x^3}\right), &
\\
&&& v_2=\epsilon\left(x^3\vfld{x^1}+\vfld{x^2}+H\vfld{x^3}\right),
\qquad
G(x)=g_0, &
\\
&&& \text{where}\quad
H=\left(\big(-c_3\big(x^3\big)^2+c_4\big)\tanh\big(c_3x^1\big)+F_{20}\big(x^2\big)\sqrt{c_3\big(x^3\big)^2-c_4}\right). &
\end{alignat*}
\end{Theorem}

\section{Optimal control problem for homogeneous metrics}

\subsection{Two states and one control}

In this section we use Pontryagin's maximum principle to compute optimal trajectories for each of the homogeneous
metrics of Theorem~\ref{2-states-1-control-homog-examples-thm}.

{\bf Case 1.1.} This point-af\/f\/ine distribution corresponds to the control system
\begin{gather}
\label{case-1-1-control-sys}
\dot{x}^{1}=1,
\qquad
\dot{x}^{2}=u,
\end{gather}
with cost functional
\begin{gather*}
Q(\dot{x})=\frac{1}{2}e^{2c_1x^1}u^2.
\end{gather*}
Consider the problem of computing optimal trajectories for~\eqref{case-1-1-control-sys}.
The Hamiltonian for the energy functional~\eqref{cost-functional-eq} is
\begin{gather*}
{\mathcal H}=p_{1}\dot{x}^{1}+p_{2}\dot{x}^{2}-Q(\dot{x})=p_{1}+p_{2}u-\frac{1}{2}e^{2c_1x^{1}}u^{2}.
\end{gather*}
By Pontryagin's maximum principle, a~necessary condition for optimal trajectories is that the control function $u(t)$
is chosen so as to maximize ${\mathcal H}$.
Since \(u\) is unrestricted and $\frac{1}{2}e^{2c_1x^{1}}>0$, \(\max_{u} {\mathcal H}\) occurs when
\begin{gather*}
0=\frac{\partial{\mathcal H}}{\partial u}=p_{2}-e^{2c_1x^{1}}u,
\end{gather*}
that is, when
\begin{gather*}
u=p_{2}e^{-2c_1x^{1}}.
\end{gather*}
So along an optimal trajectory, we have
\begin{gather*}
{\mathcal H}=p_{1}+(p_{2})^{2}e^{-2c_1x^{1}}-\frac{1}{2}(p_{2})^{2}e^{-2c_1x^{1}}=
p_{1}+\frac{1}{2}(p_{2})^{2}e^{-2c_1x^{1}}.
\end{gather*}
Moreover, \({\mathcal H}\) is constant along trajectories, and so we have
\begin{gather*}
p_{1}+\frac{1}{2}(p_{2})^{2}e^{-2c_1x^{1}}=k.
\end{gather*}
Hamilton's equations
\begin{displaymath}
\dot{x} = \frac{\partial {\mathcal H}}{\partial p},
\qquad
\dot{p} = - \frac{\partial {\mathcal H}}{\partial x}
\end{displaymath}
take the form
\begin{alignat}{3}
& \dot{x}^{1}=1,
\qquad &&
\dot{p}_{1}=c_1(p_2)^2e^{-2c_1x^1}, &\nonumber\\
&
\dot{x}^{2}=p_2e^{-2c_1x^1},
\qquad &&
\dot{p}_{2}=0. & \label{case-1-1-Ham-eqs-3}
\end{alignat}
The equation for $\dot{p_2}$ in~\eqref{case-1-1-Ham-eqs-3} implies that $p_2$ is constant; say, $p_2=c_2$.
Then optimal trajectories are solutions of the system
\begin{gather*}
\dot{x}^{1}=1,
\qquad
\dot{x}^{2}=c_2e^{-2c_1x^1}.
\end{gather*}
This system can be integrated explicitly:
\begin{itemize}\itemsep=0pt
\item If $c_1=0$, then the solutions are
\begin{gather*}
x^{1}=t,
\qquad
x^{2}=c_2t+c_3.
\end{gather*}
These solutions correspond to the family of curves
\begin{gather*}
x^2=c_2x^1+c_3
\end{gather*}
in the $\big(x^1,x^2\big)$-plane.
Thus, the set of critical curves consists of all non-vertical straight lines in the $\big(x^1,x^2\big)$ plane, oriented in the
direction of increasing $x^1$.

\item If $c_1\neq0$, then the solutions are
\begin{gather*}
x^{1}=t,
\qquad
x^{2}=-\frac{1}{2c_1}c_2e^{-2c_1t}.
\end{gather*}
These solutions correspond to the family of curves
\begin{gather*}
x^2=-\frac{1}{2c_1}c_2e^{-2c_1x^1}
\end{gather*}
in the $\big(x^1,x^2\big)$-plane.
Thus, the set of critical curves consists of a~family of exponential curves in the $\big(x^1,x^2\big)$ plane, oriented in the
direction of increasing $x^1$.
\end{itemize}

{\bf Case 1.2.} This point-af\/f\/ine distribution corresponds to the control system
\begin{gather*}
\dot{x}^{1}=x^{2},
\qquad
\dot{x}^{2}=x^{2}j_0+x^{2}u,
\end{gather*}
with cost functional
\begin{gather*}
Q(\dot{x})=\frac{1}{2}g_0u^2.
\end{gather*}

Pontryagin's maximum principle leads to the Hamiltonian
\begin{gather*}
{\mathcal H}=p_{1}x^{2}+p_{2}x^{2}j_0+\frac{1}{2g_0}\bigl(p_{2}x^{2}\bigr)^{2}
\end{gather*}
along an optimal trajectory, and Hamilton's equations take the form
\begin{alignat}{3}
& \dot{x}^{1}=x^{2}, \qquad && \dot{p}_{1}=0, &\nonumber\\
&
\dot{x}^{2}=x^{2}j_{0}+\frac{p_{2}\bigl(x^{2}\bigr)^{2}}{g_0},
\qquad &&
\dot{p}_{2}=-p_{1}-p_{2}j_{0}-\frac{(p_{2})^{2}x^{2}}{g_0}. &\label{case-1-2-Ham-eqs-2}
\end{alignat}

It is straightforward to show that the three functions
\begin{gather*}
I_1={\mathcal H}=p_{1}x^{2}+p_{2}x^{2}j_{0}+\frac{1}{2g_{0}}\big(p_{2}x^{2}\big)^{2},
\qquad
I_2=p_{1},
\qquad
I_3=p_{1}x^{1}+p_{2}x^{2}
\end{gather*}
are f\/irst integrals for this system.
This observation alone would in principle allow us to construct unparametrized solution curves for the system.
But in fact, we can solve this system fully, as follows.

The equation for $\dot{p}_1$ in~\eqref{case-1-2-Ham-eqs-2} implies that $p_1$ is constant; say, \(p_{1} = c_{1}\).
Now it is straightforward to show that
\begin{gather}
\frac{d}{dt}\bigl(p_{2}x^{2}\bigr)+c_{1}x^{2}=0.
\label{case-1-2-intermed-eq-1}
\end{gather}
If \(c_{1}=0\), then~\eqref{case-1-2-intermed-eq-1} implies that \(p_{2}x^{2}\) is equal to a~constant $k_2$, and so
\begin{gather*}
\dot{x}^{2}=x^{2}\left(j_{0}+\frac{k_{2}}{g_{0}}\right)=c_{2}x^{2}.
\end{gather*}
There are two subcases, depending on the value of $c_2$.
\begin{itemize}\itemsep=0pt
\item If \(c_2 = 0\), then
$
x^2=c_3$,
and since \(\dot{x}^{1} = x^{2}\), we have
$
x^1=c_3t+c_4$.
These solutions correspond to the family of curves
$
x^2=c_3
$
in the $\big(x^1,x^2\big)$-plane.
These curves are all horizontal lines, oriented in the direction of increasing $x^1$ when $x^2>0$ and decreasing $x^1$
when $x^2<0$.

\item If $c_2\neq0$, then
$x^{2}=c_3e^{c_{2}t}$,
and since \(\dot{x}^{1} = x^{2}\), we have
$x^{1}=\frac{c_3}{c_{2}}e^{c_{2}t}+c_{4}$.
These solutions correspond to the family of curves
$x^2=c_2\big(x^1-c_4\big)$
in the $\big(x^1,x^2\big)$-plane.
These curves are all non-vertical, non-horizontal lines, oriented in the direction of increasing $x^1$ when $x^2>0$ and
decreasing $x^1$ when $x^2<0$.
\end{itemize}
On the other hand, if \(c_{1} \ne 0\), then it is straightforward to show that
\begin{gather*}
\frac{d^2}{dt^2}\big(p_{2}x^{2}\big)=\frac{d}{dt}\big(p_{2}x^{2}\big)\biggl(j_{0}+\frac{p_{2}x^{2}}{g_{0}}\biggr).
\end{gather*}
Integrating this equation once gives
\begin{gather}
\label{case-1-2-intermed-eq-3}
\frac{d}{dt}\big(p_{2}x^{2}\big)=j_{0}\big(p_{2}x^{2}\big)+\frac{\big(p_{2}x^{2}\big)^{2}}{2g_{0}}+c_2.
\end{gather}
There are three subcases, depending on the value of $k=g_0(j_0^2g_0-2c_2)$.
\begin{itemize}\itemsep=0pt
\item If \(k = 0\), then the solution to~\eqref{case-1-2-intermed-eq-3} is
\[
p_{2} x^{2} = -\frac{g_{0}(2 + j_0(t + c_3))}{t + c_3},
\]
and from equation~\eqref{case-1-2-intermed-eq-1},
\[
x^{2} = -\frac{1}{c_{1}} \frac{d}{dt}\big(p_{2}x^{2}\big) = - \frac{2g_{0}}{c_{1}(t + c_3)^{2}}.
\]
Then since \(\dot{x}^{1} = x^{2} = - \frac{1}{c_{1}} \frac{d}{dt}\big(p_{2}x^{2}\big)\), we have
\begin{gather*}
x^{1}=-\frac{1}{c_{1}}\big(p_{2}x^{2}\big)+c_{4}=\frac{g_{0}(2+j_0(t+c_3))}{c_1(t+c_3)}+c_{4}.
\end{gather*}
These solutions correspond to the family of curves
\begin{gather*}
x^2=-\frac{1}{2c_1g_0}\left(c_1x^1-(j_0g_0+c_1c_4)\right)^2
\end{gather*}
in the $\big(x^1,x^2\big)$-plane.
These curves are all parabolas with vertex lying on the $x^1$-axis.
Since we must have $x^2\neq0$, the set of critical curves consists of all branches of parabolas with vertex on the
$x^2$-axis, oriented in the direction of increasing $x^1$ when $x^2>0$ and decreasing $x^1$ when $x^2<0$.

\item If $k>0$, then the solution to~\eqref{case-1-2-intermed-eq-3} is
\begin{gather*}
p_2x^2=-\sqrt{k}\tanh\left(\frac{\sqrt{k}}{2g_0}(t+c_3)\right)-j_0g_0,
\end{gather*}
and from equation~\eqref{case-1-2-intermed-eq-1},
\begin{gather*}
x^{2}=-\frac{1}{c_{1}}\frac{d}{dt}\big(p_{2}x^{2}\big)=\frac{k}{2c_1g_0}\text{sech}^2\left(\frac{\sqrt{k}}{2g_0}(t+c_3)\right).
\end{gather*}
Then since \(\dot{x}^{1} = x^{2} = - \frac{1}{c_{1}} \frac{d}{dt}\big(p_{2}x^{2}\big)\), we have
\begin{gather*}
x^{1}=-\frac{1}{c_{1}}\big(p_{2}x^{2}\big)+c_{4}
=\frac{1}{c_{1}}\left(\sqrt{k}\tanh\left(\frac{\sqrt{k}}{2g_0}(t+c_3)\right)+j_0g_0\right)+c_{4}.
\end{gather*}
These solutions correspond to the family of curves
\begin{gather*}
x^2=-\frac{1}{2c_1g_0}\big[\big(c_1x^1-(j_0g_0+c_1c_4)\big)^2-k\big]
\end{gather*}
in the $\big(x^1,x^2\big)$-plane.
These curves are all parabolas opening towards the $x^1$-axis.
Thus the set of critical curves consists of parabolic arcs opening towards the $x^1$-axis, approaching the axis as
$t\to\pm\infty$, and oriented in the direction of increasing $x^1$ when $x^2>0$ and decreasing $x^1$ when $x^2<0$.

\item If $k<0$, then the solution to~\eqref{case-1-2-intermed-eq-3} is
\begin{gather*}
p_2x^2=\sqrt{-k}\tan\left(\frac{\sqrt{-k}}{2g_0}(t+c_3)\right)-j_0g_0,
\end{gather*}
and from equation~\eqref{case-1-2-intermed-eq-1},
\begin{gather*}
x^{2}=-\frac{1}{c_{1}}\frac{d}{dt}\big(p_{2}x^{2}\big)=\frac{k}{2c_1g_0}\text{sec}^2\left(\frac{\sqrt{-k}}{2g_0}(t+c_3)\right).
\end{gather*}
Then since \(\dot{x}^{1} = x^{2} = - \frac{1}{c_{1}} \frac{d}{dt}\big(p_{2}x^{2}\big)\), we have
\begin{gather*}
x^{1}=-\frac{1}{c_{1}}\big(p_{2}x^{2}\big)+c_{4}
=-\frac{1}{c_{1}}\left(\sqrt{-k}\tan\left(\frac{\sqrt{-k}}{2g_0}(t+c_3)\right)-j_0g_0\right)+c_{4}.
\end{gather*}
These solutions correspond to the family of curves
\begin{gather*}
x^2=-\frac{1}{2c_1g_0}\big[\big(c_1x^1-(j_0g_0+c_1c_4)\big)^2-k\big]
\end{gather*}
in the $\big(x^1,x^2\big)$-plane.
These curves are all parabolas opening away from the $x^1$-axis.
Thus the set of critical curves consists of parabolic arcs opening away from the $x^1$-axis, becoming unbounded in
f\/inite time, and oriented in the direction of increasing $x^1$ when $x^2>0$ and decreasing $x^1$ when $x^2<0$.
\end{itemize}

\subsection{Three states and one control}

In this section we use Pontryagin's maximum principle to compute optimal trajectories for each of the homogeneous
metrics of Theorem~\ref{3-states-1-control-homog-examples-thm}.

{\bf Case 2.1.1.} This point-af\/f\/ine distribution corresponds to the control system
\begin{gather*}
\dot{x}^{1}=1,
\qquad
\dot{x}^{2}=x^3,
\qquad
\dot{x}^3=c_2x^2+c_3x^3+u,
\end{gather*}
with cost functional
\begin{gather*}
Q(\dot{x})=\frac{1}{2}u^2.
\end{gather*}

\begin{figure}[t] \centering
\includegraphics[width=2.5in]{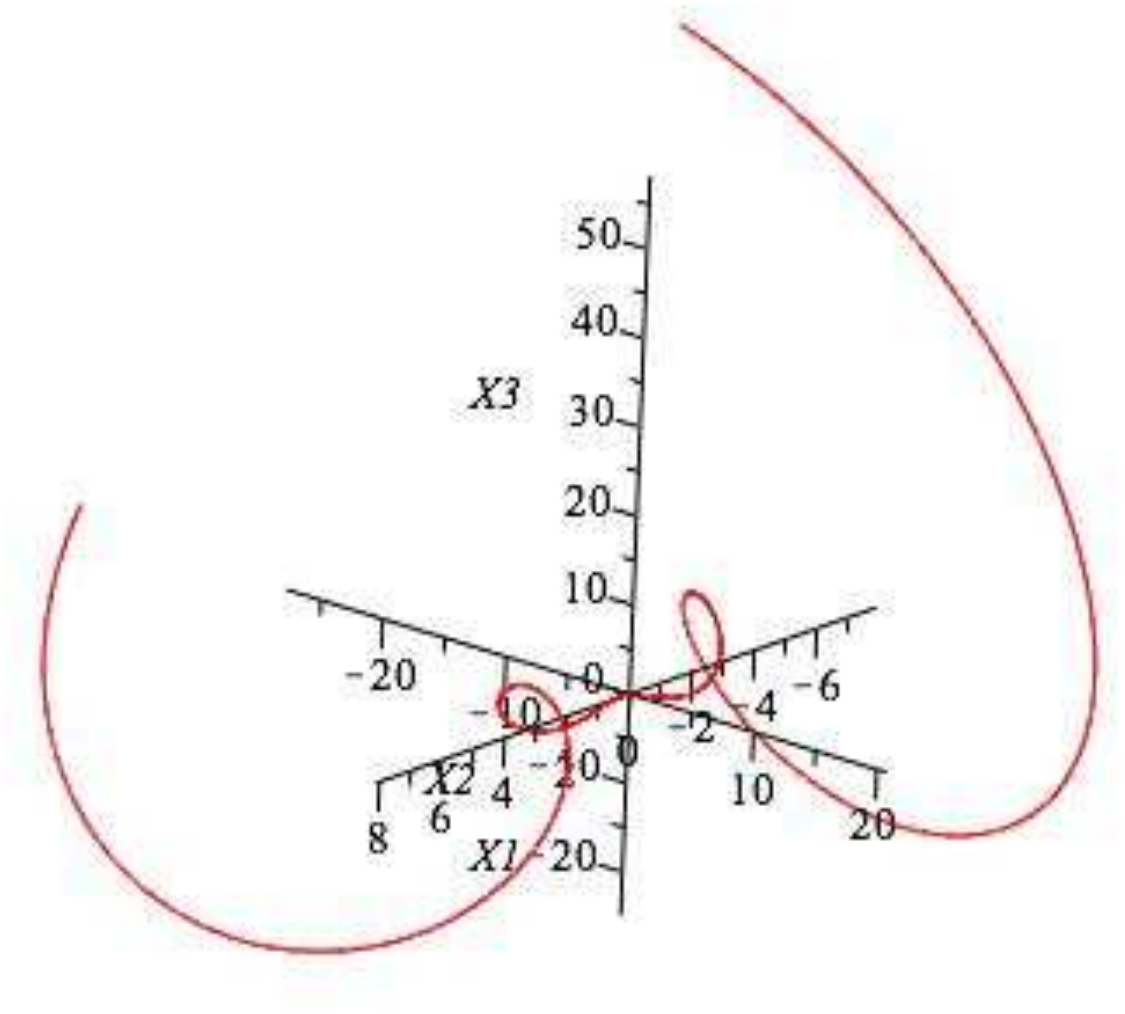}
 \caption{}\label{fig:case-2-1-1}
\end{figure}

The Hamiltonian for the energy functional~\eqref{cost-functional-eq} is
\begin{gather*}
{\mathcal H}=p_{1}\dot{x}^{1}+p_{2}\dot{x}^{2}+p_{3}\dot{x}^{3}-Q(\dot{x})
=p_{1}+p_{2}x^{3}+p_3\big(c_2x^2+c_3x^3+u\big)-\frac{1}{2}u^{2}.
\end{gather*}
Pontryagin's maximum principle leads to the Hamiltonian
\begin{gather*}
{\mathcal H}=p_{1}+p_{2}x^{3}+p_3\big(c_2x^2+c_3x^3\big)+\frac{1}{2}(p_3)^2
\end{gather*}
along an optimal trajectory, and Hamilton's equations take the form
\begin{alignat}{3}
& \dot{x}^{1}=1,
\qquad &&
\dot{p}_{1}=0, &
\nonumber\\
& \dot{x}^{2}=x^3,
\qquad &&
\dot{p}_{2}=-c_2p_3, &
\nonumber
\\
& \dot{x}^3=c_2x^2+c_3x^3+p_3,
\qquad &&
\dot{p}_3=-p_2-c_3p_3. &
\label{case-2-1-1-Ham-eqs}
\end{alignat}
The equations for $\dot{p}_2$ and $\dot{p}_3$ in~\eqref{case-2-1-1-Ham-eqs} can be written as
\begin{gather*}
\ddot{p}_2+c_3\dot{p}_2-c_2p_2=0,
\end{gather*}
and the function $p_3=-\frac{1}{c_2}\dot{p}_2$ satisf\/ies this same ODE.
Then the equations for $\dot{x}^2$ and $\dot{x}^3$ can be written as
\begin{gather*}
\ddot{x}^2-c_3\dot{x}^2-c_2x^2=p_3(t),
\end{gather*}
where $p_3(t)$ is an arbitrary solution of the ODE
\begin{gather*}
\ddot{p}_3+c_3\dot{p}_3-c_2p_3=0.
\end{gather*}
Therefore, $x^2(t)$ is an arbitrary solution of the 4th-order ODE
\begin{gather*}
\left(\frac{d^2}{dt^2}+c_3\frac{d}{dt}-c_2\right)\left(\frac{d^2}{dt^2}-c_3\frac{d}{dt}-c_2\right)x^2(t)=0,
\end{gather*}
and for any such $x^2(t)$, we have
\begin{gather*}
x^1(t)=t+t_0,
\qquad
x^3(t)=\dot{x}^2(t).
\end{gather*}
A sample optimal trajectory is shown in Fig.~\ref{fig:case-2-1-1}.

{\bf Case 2.1.2.} This point-af\/f\/ine distribution corresponds to the control system
\begin{gather*}
\dot{x}^{1}=1,
\qquad
\dot{x}^{2}=x^3,
\qquad
\dot{x}^3=c_3x^3+u,
\end{gather*}
with cost functional
\begin{gather*}
Q(\dot{x})=\frac{1}{2\big(c_1x^3\big)^2}u^2.
\end{gather*}

\begin{figure}[t] \centering
\includegraphics[width=2.7in]{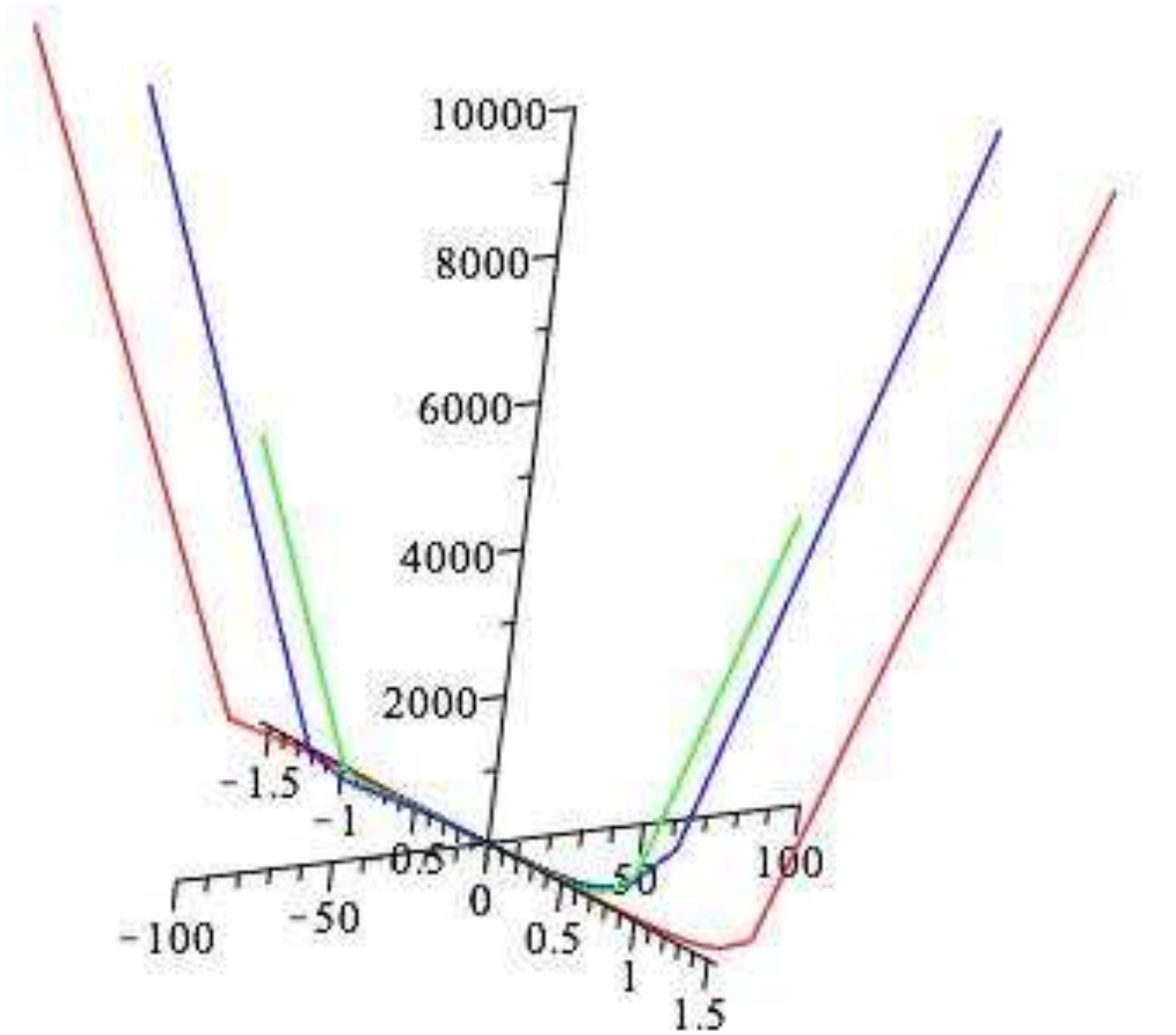} \includegraphics[width=1.9in]{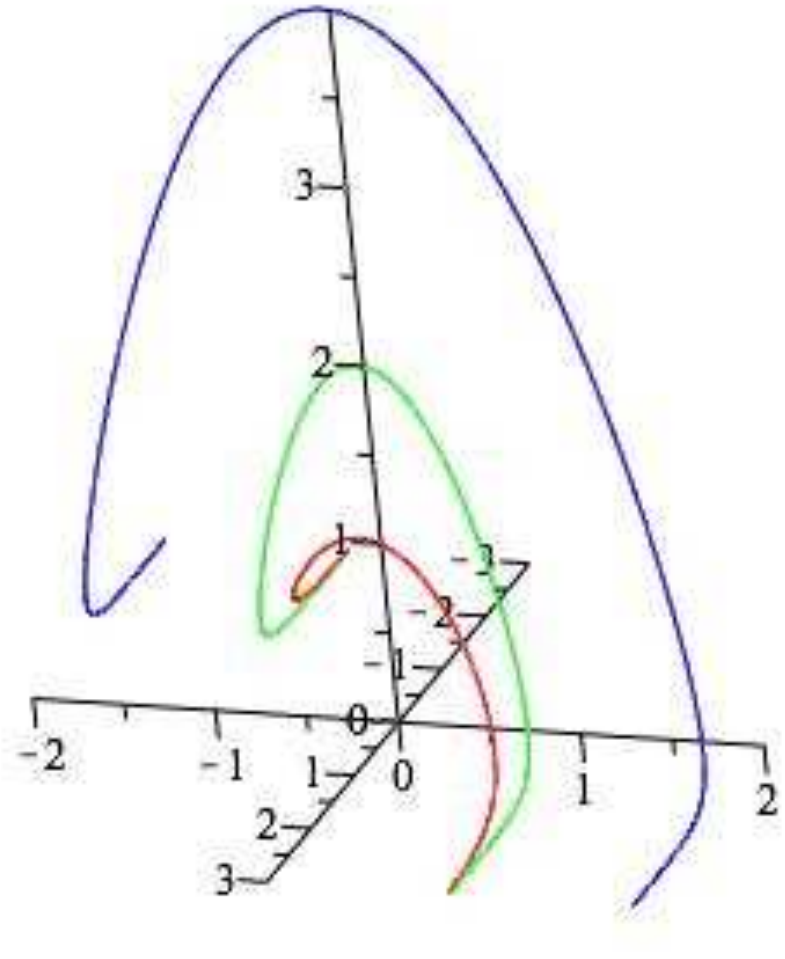} \caption{}
\label{fig:case-2-1-2-c2nonzero}
\end{figure}

Pontryagin's maximum principle leads to the Hamiltonian
\begin{gather*}
{\mathcal H}=p_{1}+p_{2}x^{3}+c_3p_3x^3+\frac{1}{2}\big(c_1x^3p_3\big)^2
\end{gather*}
along an optimal trajectory, and Hamilton's equations take the form
\begin{alignat}{3}
& \dot{x}^{1}=1,
\qquad &&
\dot{p}_{1}=0, &
\nonumber\\
& \dot{x}^{2}=x^3,
\qquad &&
\dot{p}_{2}=0,&
\nonumber
\\
& \dot{x}^3=c_3x^3+\big(c_1x^3\big)^2p_3,
\qquad &&
\dot{p}_3=-p_2-c_3p_3-(c_1p_3)^2x^3. &
\label{case-2-1-2-Ham-eqs}
\end{alignat}
The equation for $\dot{p}_2$ in~\eqref{case-2-1-2-Ham-eqs} implies that $p_2(t)$ is equal to a~constant $c_2$.
Then~\eqref{case-2-1-2-Ham-eqs} implies that
\begin{gather}\label{case-2-1-2-almost done}
\dot{\big(p_3x^3\big)}=-c_2x^3,
\qquad
\dot{x}^3=c_3x^3+c_1x^3\big(p_3x^3\big).
\end{gather}
These equations can be solved as follows:
\begin{itemize}
\itemsep=0pt \item If $c_2=0$, then the function $p_3x^3$ is constant, and so the equation for $\dot{x}^3$ becomes
\begin{gather*}
\dot{x}^3=\tilde{c}x^3
\end{gather*}
for some constant $\tilde{c}$.
If $\tilde{c}=0$, then the solution trajectories are given by
\begin{gather*}
x^1(t)=t+t_0,
\qquad
x^2(t)=at+b,
\qquad
x^3(t)=a
\end{gather*}
for some constants $a$, $b$.
Sample optimal trajectories are shown in Fig.~\ref{fig:case-2-1-2-c2zero-ctildezero}.
\begin{figure}[t] \centering
\begin{minipage}[b]{75mm}
\includegraphics[width=2.7in]{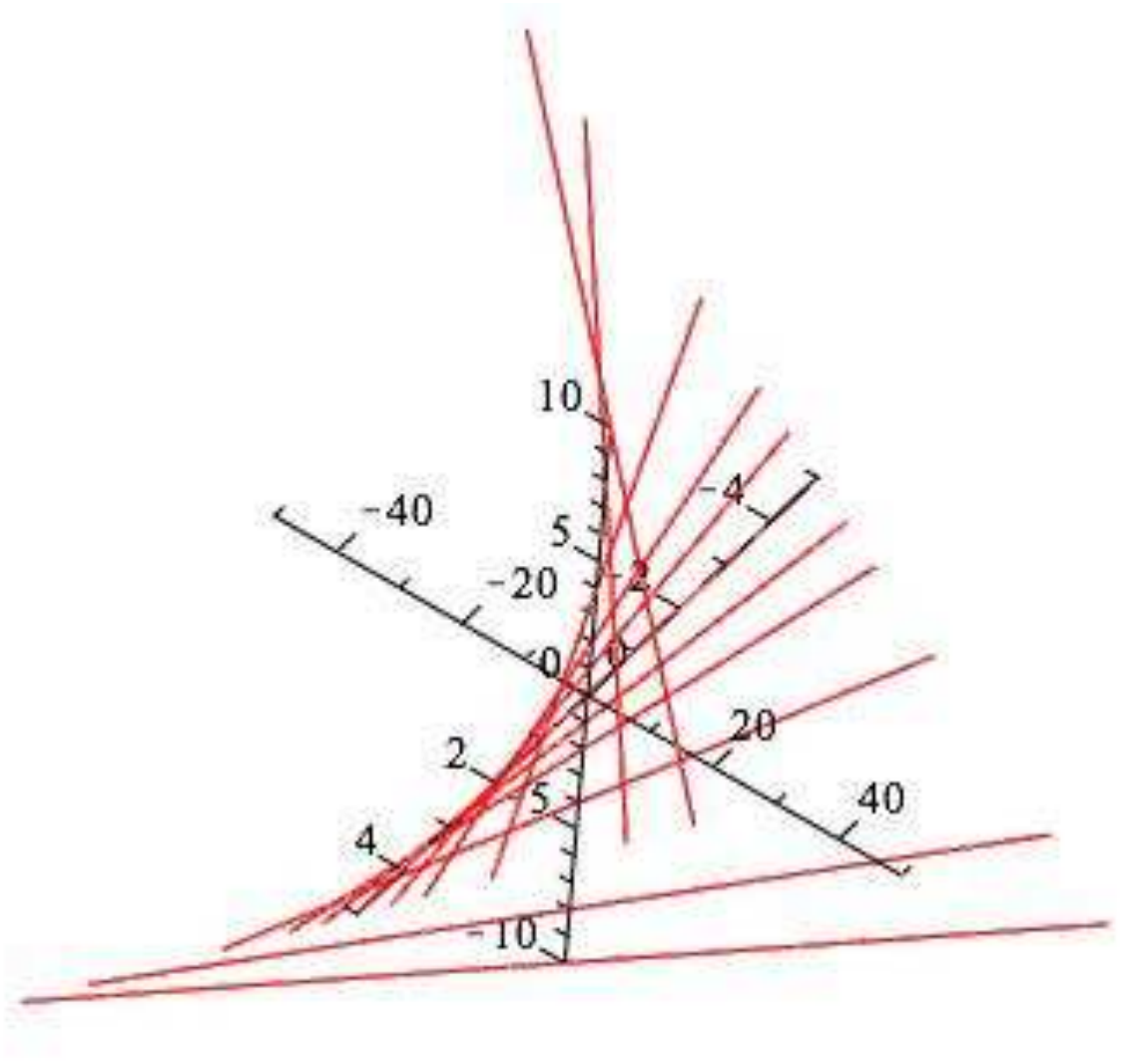} \caption{}
\label{fig:case-2-1-2-c2zero-ctildezero}
\end{minipage}
\quad
\begin{minipage}[b]{75mm}
\includegraphics[width=2.9in]{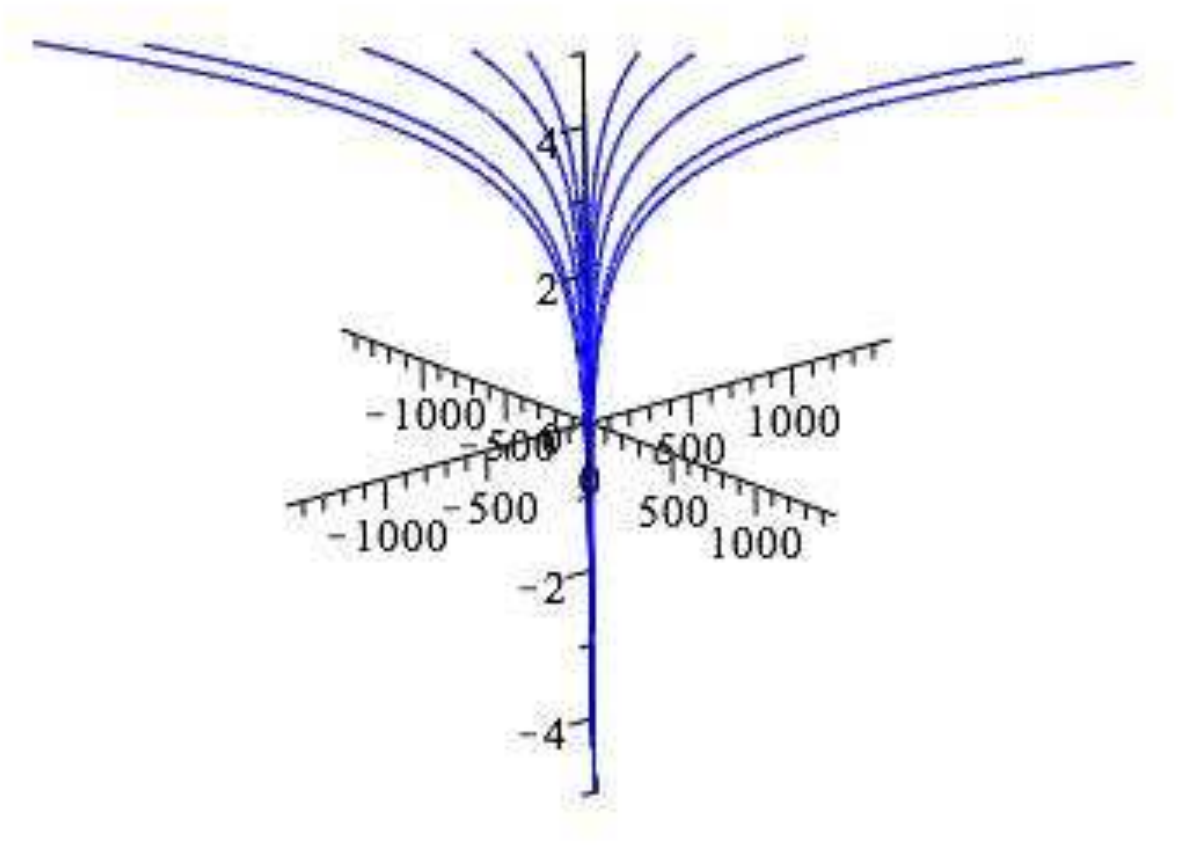} \caption{}
\label{fig:case-2-1-2-c2zero-ctildenonzero}
\end{minipage}
\end{figure}

If $\tilde{c}\neq0$, then the solution trajectories are given by
\begin{gather*}
x^1(t)=t+t_0,
\qquad
x^2(t)=\frac{a}{\tilde{c}}e^{\tilde{c}t}+b,
\qquad
x^3(t)=a e^{\tilde{c}t}
\end{gather*}
for some constants $a$, $b$.
Sample optimal trajectories are shown in Fig.~\ref{fig:case-2-1-2-c2zero-ctildenonzero}.

\item If $c_2\neq0$, then~\eqref{case-2-1-2-almost done} can be written as the 2nd-order ODE for the function
$z(t)=p_3(t)x^3(t)$:
\begin{gather*}
\ddot{z}=\big(c_3+c_1^2z\big)\dot{z}.
\end{gather*}
Integrating once yields
\begin{gather*}
\dot{z}=\frac{1}{2}(c_1z)^2+c_3z+c_4
\end{gather*}
for some constant $c_4$.
Depending on the values of the constants, the solution $z(t)$ has one of the following forms:
\begin{alignat*}{4}
&  (1) \quad &&
 z(t)  = a \tan (bt + c) + d,
\qquad &&
 c_3^2 - 2 c_1 c_4 < 0; &
\\
& (2)
\quad &&
 z(t)  = a \tanh (bt + c) + d,
\qquad &&
 c_3^2 - 2 c_1 c_4 > 0; &
\\
& (3)
\quad &&
 z(t)  = \frac{1}{at + b} + c,
\qquad &&
 c_3^2 - 2 c_1 c_4 = 0. &
\end{alignat*}
Then we have
\begin{gather*}
x^3=-\frac{1}{c_2}\dot{z}=\dot{x}^2,
\end{gather*}
and so the corresponding solution trajectories are given (with slightly modif\/ied constants) by:
\begin{gather*}
 (1)
\quad
\begin{cases}
x^1(t) = t + t_0,
\\
x^2(t) = a \tan (bt + c) + d,
\\
x^3(t) = ab \sec^2(bt + c);
\end{cases}
\qquad
 (2)
\quad
\begin{cases}
x^1(t) = t + t_0,
\\
x^2(t) = a \tanh (bt + c) + d,
\\
x^3(t) = ab \,\text{sech}^2(bt + c);
\end{cases}
\\
 (3)
\quad
\begin{cases}
x^1(t) = t + t_0,
\\
x^2(t) = \dfrac{1}{at + b} + c,
\vspace{1mm}\\
x^3(t) = -\dfrac{a}{(at + b)^2}.
\end{cases}
\end{gather*}
\end{itemize}
Sample optimal trajectories for the f\/irst two cases are shown in Fig.~\ref{fig:case-2-1-2-c2nonzero}.

{\bf Case 2.2.} This point-af\/f\/ine distribution corresponds to the control system
\begin{gather*}
\dot{x}^{1}=x^2,
\qquad
\dot{x}^{2}=x^3,
\qquad
\dot{x}^3=x^2\left(\frac{3}{2}\left(\frac{x^3}{x^2}\right)^2+c_1+u\right),
\end{gather*}
with cost functional
\begin{gather*}
Q(\dot{x})=\frac{1}{2}g_0u^2.
\end{gather*}

Pontryagin's maximum principle leads to the Hamiltonian
\begin{gather*}
{\mathcal H}
=p_{1}x^2+p_{2}x^{3}+p_3x^2\left(\frac{3}{2}\left(\frac{x^3}{x^2}\right)^2+c_1\right)+\frac{1}{2g}\big(p_3x^2\big)^2
\end{gather*}
along an optimal trajectory, and Hamilton's equations take the form
\begin{alignat}{3}
& \dot{x}^{1}=x^2,
\qquad && \dot{p}_{1}=0,&
\nonumber\\
& \dot{x}^{2}=x^3,
\qquad && \dot{p}_{2}=-p_1+\frac{3}{2}\frac{p_3\big(x^3\big)^2}{\big(x^2\big)^2}-c_1p_3-\frac{1}{g}(p_3)^2x^2,&
\nonumber\\
& \dot{x}^3=\frac{3}{2}\frac{\big(x^3\big)^2}{x^2}+c_1x^2+\frac{1}{g}p_3\big(x^2\big)^2,
\qquad &&
\dot{p}_3=-p_2-3\frac{p_3x^3}{x^2}.&
\label{case-2-2-Ham-eqs}
\end{alignat}

\begin{figure}[t] \centering \includegraphics[width=2.5in]{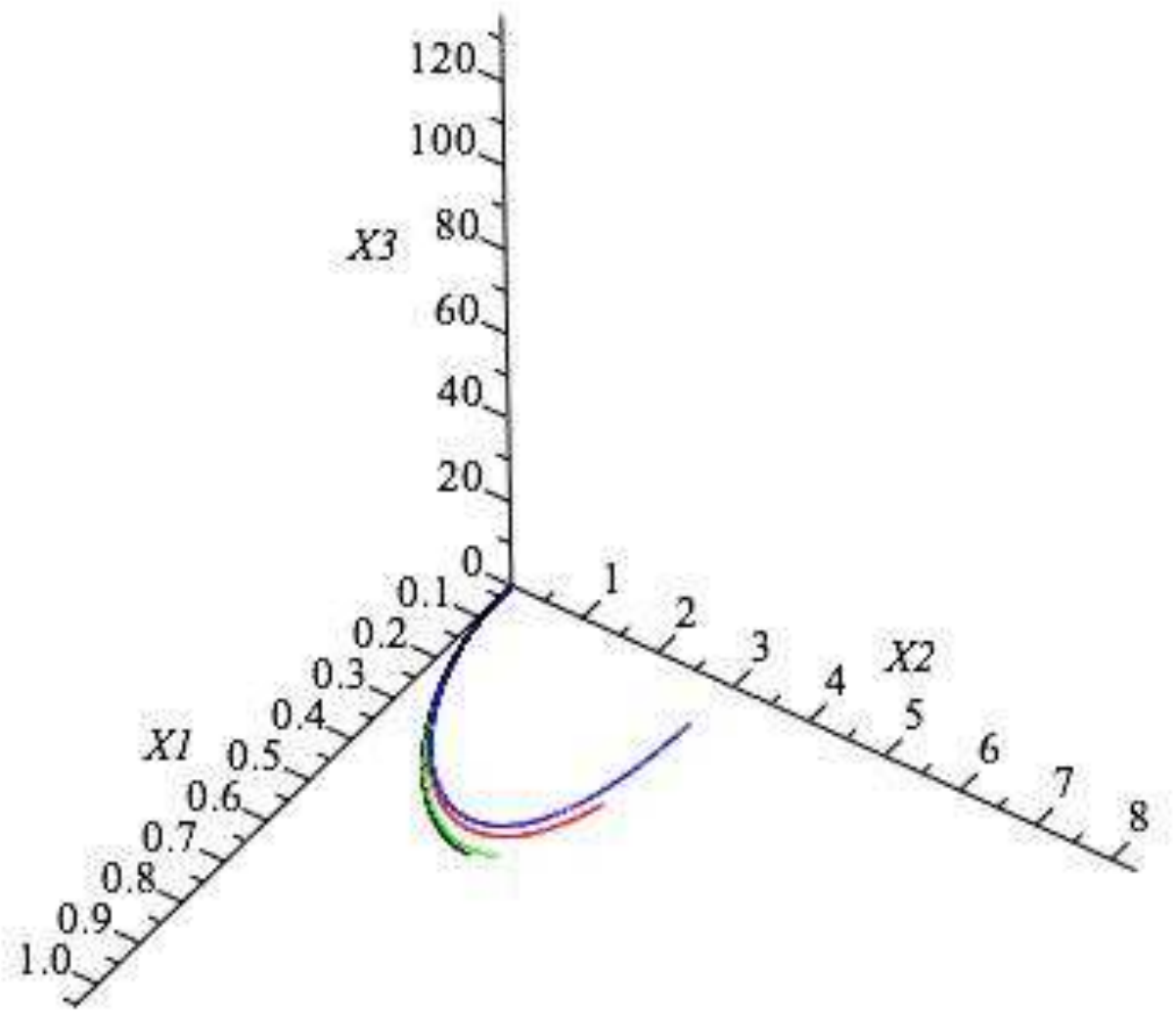} \caption{}
\label{fig:case-2-2}
\end{figure}

The system~\eqref{case-2-2-Ham-eqs} has three independent f\/irst integrals in addition to the Hamiltonian ${\mathcal
H}$ (which is automatically a~f\/irst integral): it is straightforward to show, using~\eqref{case-2-2-Ham-eqs}, that
the functions
\begin{gather*}
I_1=p_1,
\qquad
I_2=p_1x^1+p_2x^2+p_3x^3,
\qquad
I_3=p_1\big(x^1\big)^2+2p_2x^1x^2+2p_3x^1x^3+2p_3\big(x^2\big)^2
\end{gather*}
are f\/irst integrals for this system.
We can use these conserved quantities to reduce the system~\eqref{case-2-2-Ham-eqs}, as follows: on any solution curve
of~\eqref{case-2-2-Ham-eqs}, we have
\begin{gather*}
I_1=k_1,
\qquad
I_2=k_2,
\qquad
I_3=k_3
\end{gather*}
for some constants $k_1$, $k_2$, $k_3$.
These equations can be solved for $p_1$, $p_2$, $p_3$ to obtain
\begin{gather*}
 p_1=k_1,
\\
p_2=
k_1\left(-\frac{x^1}{x^2}-\frac{(x^1)^2x^3}{2(x^2)^3}\right)
+k_2\left(\frac{1}{x^2}+\frac{x^1x^3}{(x^2)^3}\right)+k_3\left(-\frac{x^3}{2(x^2)^3}\right),
\\
 p_3=
k_1\left(\frac{(x^1)^2}{2(x^2)^2}\right)+k_2\left(-\frac{x^1}{(x^2)^2}\right)
+k_3\left(\frac{1}{2(x^2)^2}\right).
\end{gather*}
These equations can be substituted into~\eqref{case-2-2-Ham-eqs} to obtain a~closed, f\/irst-order ODE system for the
functions $x^1$, $x^2$, $x^3$, depending on the parameters $k_1$, $k_2$, $k_3$; moreover, making the same substitution
in the Hamiltonian ${\mathcal H}$ yields a~conserved quantity for this system.
(The precise expressions for the system and the conserved quantity are complicated and unenlightening, so we will not
write them out explicitly here.) The resulting ODE system cannot be solved analytically, but numerical integration
yields sample trajectories as shown in Fig.~\ref{fig:case-2-2}.

{\bf Case 2.3.1.} This point-af\/f\/ine distribution corresponds to the control system
\begin{gather*}
\dot{x}^{1}=1+x^3(c_1+u),
\qquad
\dot{x}^{2}=c_1+u,
\qquad
\dot{x}^3=\epsilon\big(x^1+c_2x^3\big)(c_1+u),
\end{gather*}
with cost functional
\begin{gather*}
Q(\dot{x})=\frac{1}{2}g_0u^2.
\end{gather*}

Pontryagin's maximum principle leads to the Hamiltonian
\begin{gather*}
\begin{split}
& {\mathcal H}=p_1+p_3x^1+\frac{c_1}{\sqrt{c_2}}\left(p_1x^3+p_2+c_3p_3x^3\right)
\\
& \phantom{{\mathcal H}=}
{}+\big(p_1x^3+p_2+c_3p_3x^3\big)\left(c_2\big(p_1x^3+p_2+p_3x^1+c_3p_3x^3\big)\right)
\end{split}
\end{gather*}
along an optimal trajectory, and Hamilton's equations take the form
\begin{gather}
\dot{x}^{1}=1+x^3\left(\frac{c_1}{\sqrt{c_2}}+c_2\big(p_1x^3+p_2+p_3x^1+c_3p_3x^3\big)\right),
\nonumber
\\
\dot{x}^{2}=\frac{c_1}{\sqrt{c_2}}+c_2\big(p_1x^3+p_2+p_3x^1+c_3p_3x^3\big),
\nonumber
\\
\dot{x}^3=\big(x^1+c_3x^3\big)\left(\frac{c_1}{\sqrt{c_2}}+c_2\big(p_1x^3+p_2+p_3x^1+c_3p_3x^3\big)\right),
\nonumber
\\
\dot{p}_{1}=-p_3\left(\frac{c_1}{\sqrt{c_2}}+c_2\big(p_1x^3+p_2+p_3x^1+c_3p_3x^3\big)\right),
\nonumber
\\
\dot{p}_{2}=0,
\nonumber
\\
\dot{p}_3=-\Biggl(\left(p_1+c_3p_3\right)\left(\frac{c_1}{\sqrt{c_2}}+c_2\big(p_1x^3+p_2+p_3x^1+c_3p_3x^3\big)\right)
\nonumber
\\
\phantom{\dot{p}_3=}
{}+\left((p_1+c_3p_3)x^3+p_2\right)\left(c_2p_1+c_2c_3p_3\right)\Biggr).
\label{case-2-3-Ham-eqs}
\end{gather}

\begin{figure}[t] \centering \includegraphics[width=1.9in]{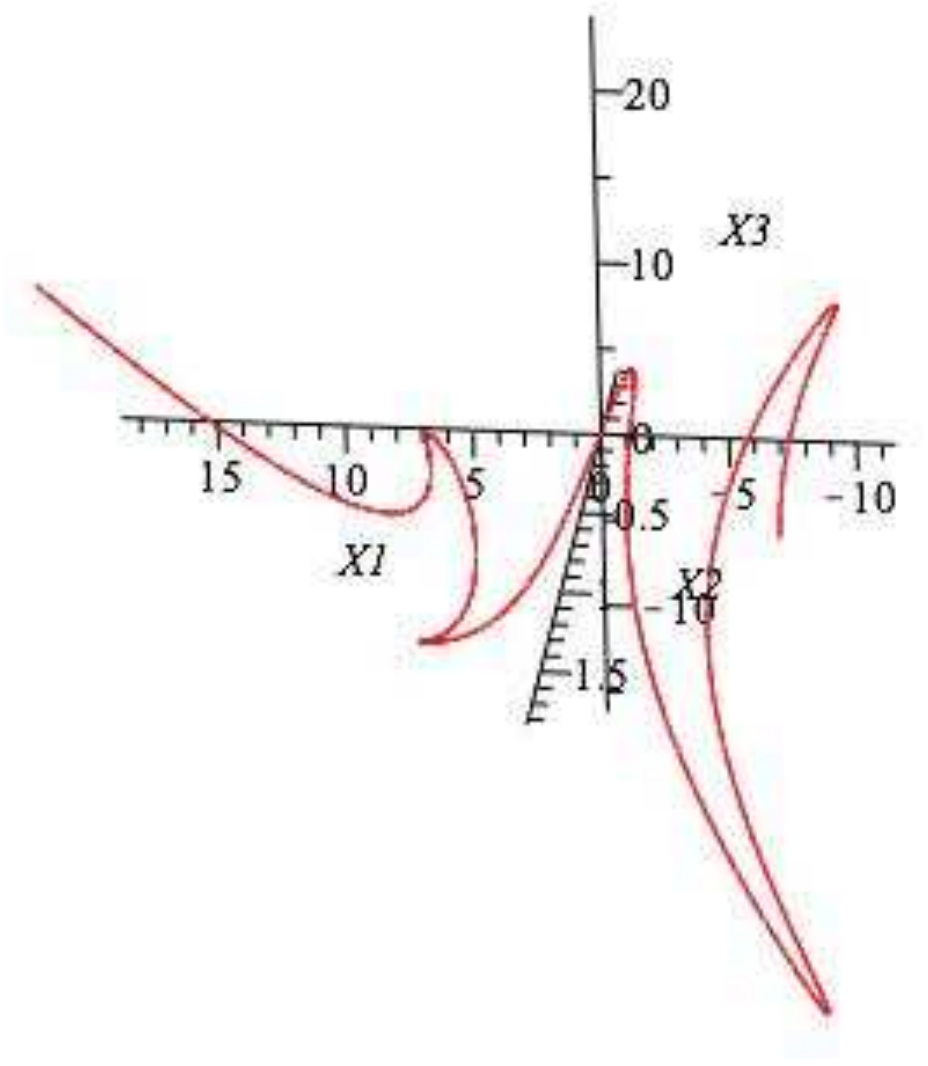}\qquad \includegraphics[width=2.2in]{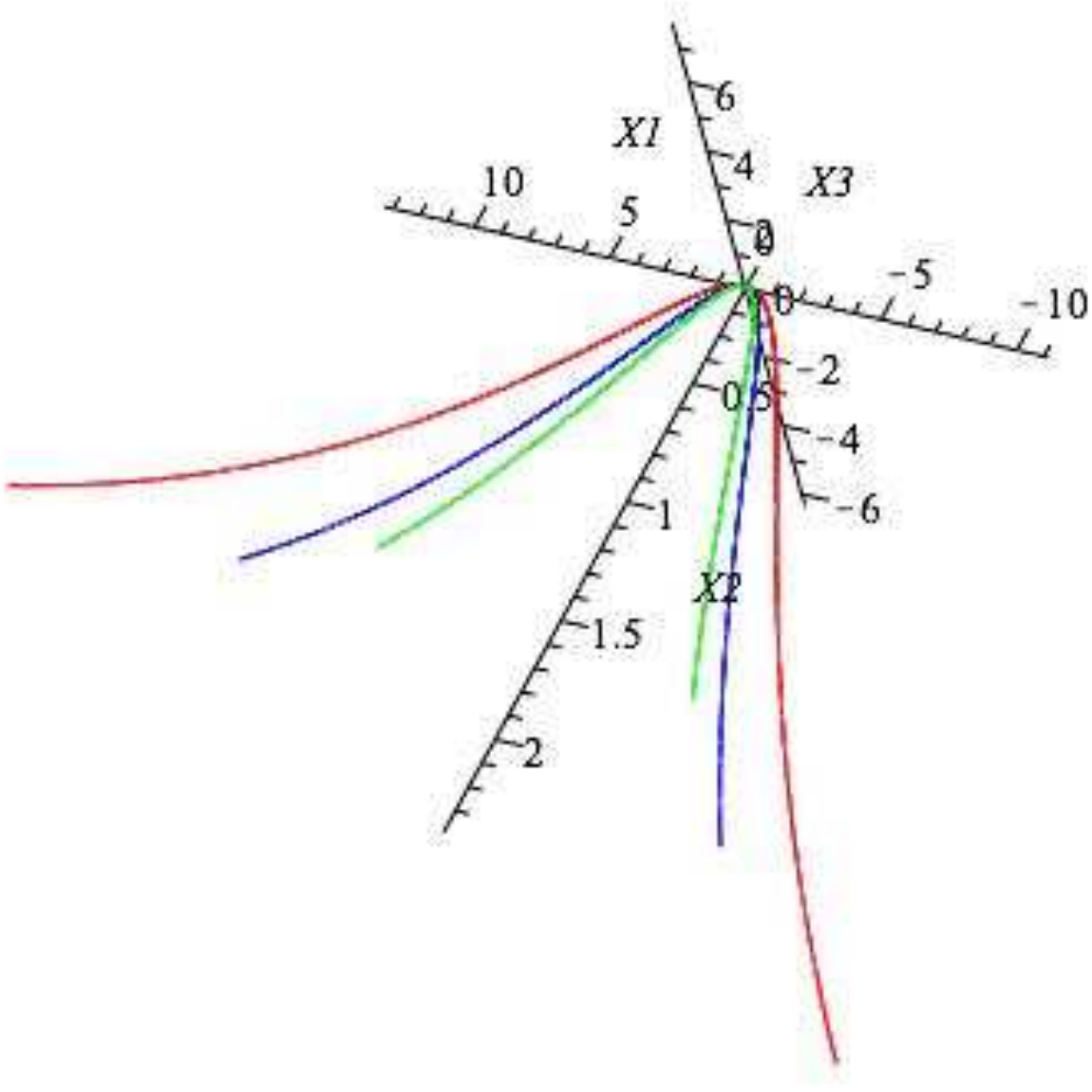} \caption{}
\label{fig:case-2-3}
\end{figure}

The system~\eqref{case-2-3-Ham-eqs} has three independent f\/irst integrals in addition to the Hamiltonian ${\mathcal
H}$: it is straightforward to show, using~\eqref{case-2-3-Ham-eqs}, that the functions
\begin{gather*}
I_1=(p_1+r_1p_3)e^{r_1x^2},
\qquad
I_2=(p_1+r_2p_3)e^{r_2x^2},
\qquad
I_3=p_2,
\end{gather*}
where $r_1$, $r_2$ are as in~\eqref{case-2-3-define-r's}, are f\/irst integrals for this system.
A similar process to that described in the previous case leads to a~closed, f\/irst-order system of ODEs for the
functions~$x^1$,~$x^2$,~$x^3$; numerical integration of this system yields sample trajectories (with $c_3>0$ and $c_3<0$) as
shown in Fig.~\ref{fig:case-2-3}.

{\bf Case 2.3.2.} Due to the complexity of the computations, we will restrict our attention to the simplest
case, where $\epsilon=1$, $c_1=c_4=0$, and $F_{20}\big(x^2\big)=0$.
(In the interest of brevity, we will omit Case~2.3.3, which is similar to this case.) With these assumptions, we have
\begin{gather*}
v_1=\vfld{x^1},
\qquad
v_2=x^3\vfld{x^1}+\vfld{x^2}+c_3\big(x^3\big)^2\tan\big(c_3x^1\big)\vfld{x^3}.
\end{gather*}
This point-af\/f\/ine distribution corresponds to the control system
\begin{gather*}
\dot{x}^{1}=1+x^3u,
\qquad
\dot{x}^{2}=u,
\qquad
\dot{x}^3=c_3\big(x^3\big)^2\tan\big(c_3x^1\big)u  ,
\end{gather*}
with cost functional
\begin{gather*}
Q(\dot{x})=\frac{1}{2}g_0u^2.
\end{gather*}

\begin{figure}[t] \centering \includegraphics[width=1.8in]{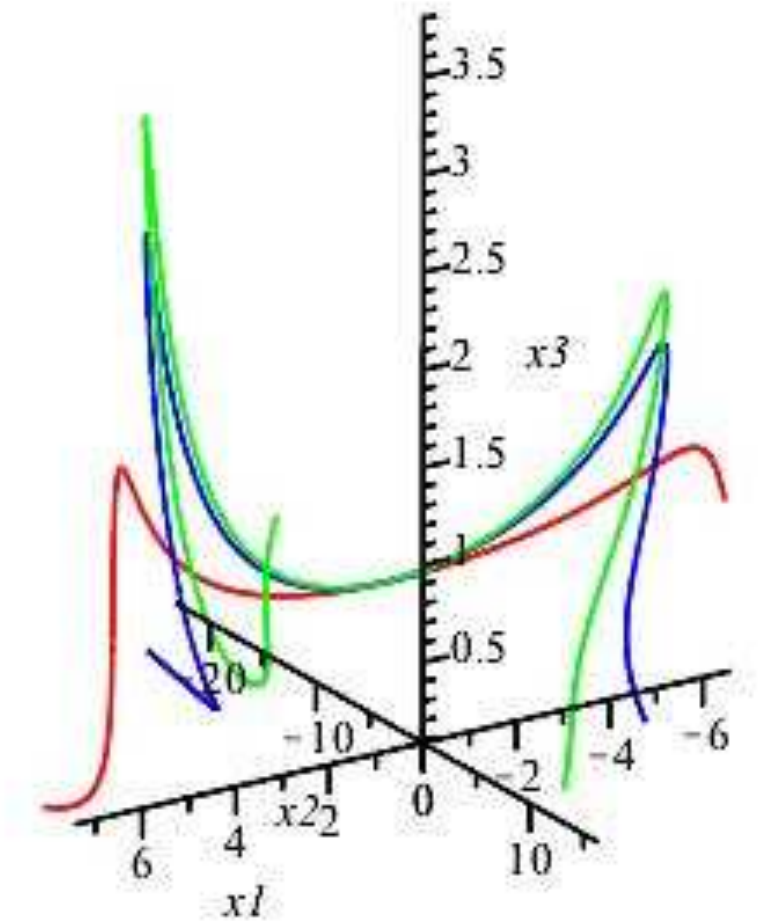} \caption{}
\label{fig:case-2-3-2}
\end{figure}

Pontryagin's maximum principle leads to the Hamiltonian
\begin{gather*}
{\mathcal H}
=\frac{1}{2g_0}\Big(c_3^2\big(\big(x^3\big)^2p_3\big)^2\tan^2\big(c_3x^1\big)
+2c_3p_3\big(x^3\big)^2\big(x^3p_1+p_2\big)\tan\big(c_3x^1\big)
\\
\phantom{{\mathcal H}=}
{}+\big(2p_1\big(x^3p_2+g_0\big)+\big(x^3p_1\big)^2+p_2^2\big)\Big)
\end{gather*}
along an optimal trajectory, and Hamilton's equations take the form
\begin{gather*}
\dot{x}^1=\frac{1}{g_0}\left(c_3\big(x^3\big)^3p_3\tan\big(c_3x^1\big)+\big(x^3\big)^2p_1+x^3p_2+g_0\right),
\\
\dot{x}^2=\frac{1}{g_0}\left(c_3\big(x^3\big)^2p_3\tan\big(c_3x^1\big)\right)+x^3p_1+p_2,
\\
\dot{x}^3=\frac{1}{g_0}c_3\big(x^3\big)^2\tan\big(c_3x^1\big)
\left(c_3\big(x^3\big)^2p_3\tan\big(c_3x^1\big)+x^3p_1+p_2\right),
\\
\dot{p}_{1}=-\frac{1}{g_0\cos^3\big(c_3x^1\big)}\big(c_3x^3\big)^2p_3
\left(c_3\big(x^3\big)^2p_3\sin\big(c_3x^1\big)+\big(x^3p_1+p_2\big)\cos\big(c_3x^1\big)\right),
\\
\dot{p}_{2}=0,
\\
\dot{p}_3=-\frac{1}{g_0}\left(2(c_3p_3)^2\big(x^3\big)^3\tan^2\big(c_3x^1\big)
+c_3x^3p_3\big(3x^3p_1+2p_2\big)\tan\big(c_3x^1\big)
+p_1\big(x^3p_1+p_2\big)\right).
\end{gather*}
This ODE system cannot be solved analytically, but numerical integration yields sample trajectories as shown in
Fig.~\ref{fig:case-2-3-2}.

\section{Conclusion}

What is perhaps most interesting about these results is how the behavior of control-af\/f\/ine systems in low
dimensions varies from that of control-linear (i.e., driftless) systems.
As we observed in~\cite{CMW09}, functional invariants appear in much lower dimension for af\/f\/ine distributions
(beginning with $n=2$, $s=1$) than for linear distributions, where there are no functional invariants in dimensions below
$n=5$, $s=2$.

With the addition of a~quadratic cost functional, we see a~similar phenomenon: for linear distributions with
a~quadratic cost functional, there are no functional invariants for any $n$ when $s=1$, since local coordinates can
always be chosen so that a~unit vector f\/ield for the cost functional is represented by the vector f\/ield
$\frac{\partial}{\partial x^1}$.
But for af\/f\/ine distributions with $s=1$, there are numerous functional invariants, and even the homogeneous
examples exhibit a~wide variety of behaviors for the optimal trajectories.

\appendix
\section{Normal forms for Case 2.3}
\label{details-app}

In this appendix, we carry out the analysis to identify examples of normal forms for homogeneous point-af\/f\/ine
structures in Case 2.3.

First consider local coordinate transformations which preserve the expressions~\eqref{case-2-3-coframing}.
Let $\big(x^1,x^2,x^3\big)$ and $(\tilde{x}^1,\tilde{x}^2,\tilde{x}^3)$ be two local coordinate systems with respect to which
the coframing $(\eta^1,\eta^2,\eta^3)$ takes the form~\eqref{case-2-3-coframing}.
Then we must have
\begin{gather}
\eta^1=dx^1-x^3\,dx^2=d\tilde{x}^1-\tilde{x}^3\,d\tilde{x}^2.
\label{case-2-3-eta1-in-coords}
\end{gather}
Taking the exterior derivative of~\eqref{case-2-3-eta1-in-coords} yields
\begin{gather}
d\eta^1=dx^2\wedge dx^3=d\tilde{x}^2\wedge d\tilde{x}^3.
\label{case-2-3-deta1-in-coords}
\end{gather}
In particular,
\begin{gather*}
\text{span}\big(dx^2,dx^3\big)=\;\text{span}\big(d\tilde{x}^2,d\tilde{x}^3\big).
\end{gather*}
Therefore we must have
\begin{gather}
x^2=\bar{\phi}\big(\tilde{x}^2,\tilde{x}^3\big),
\qquad
x^3=\bar{\psi}\big(\tilde{x}^2,\tilde{x}^3\big)
\label{case-2-3-bad-coord-trans}
\end{gather}
for some functions $\bar{\phi}\big(\tilde{x}^2,\tilde{x}^3\big)$, $\bar{\psi}\big(\tilde{x}^2,\tilde{x}^3\big)$.
Equation~\eqref{case-2-3-deta1-in-coords} then implies that the functions~$\bar{\phi}$,~$\bar{\psi}$ satisfy the PDE
\begin{gather}
\bar{\phi}_{\tilde{x}^2}\bar{\psi}_{\tilde{x}^3}-\bar{\phi}_{\tilde{x}^3}\bar{\psi}_{\tilde{x}^2}=1.
\label{case-2-3-bad-transformation-pde}
\end{gather}

Unfortunately, equation~\eqref{case-2-3-bad-transformation-pde} cannot be solved explicitly in terms of arbitrary
functions of~$\tilde{x}^2$,~$\tilde{x}^3$.
However, it {\em can} be solved implicitly with a~slightly dif\/ferent setup.
Instead of~\eqref{case-2-3-bad-coord-trans}, suppose that we def\/ine our coordinate transformation by
\begin{gather*}
\tilde{x}^2=\phi\big(x^2,\tilde{x}^3\big),
\qquad
x^3=\psi\big(x^2,\tilde{x}^3\big).
\end{gather*}
Then equation~\eqref{case-2-3-deta1-in-coords} is equivalent to the condition
\begin{gather*}
\phi_{x^2}=\psi_{\tilde{x}^3}.
\end{gather*}
(In addition, both terms in this equation must be nonzero.) This is equivalent to the condition that there exists
a~function $\Phi\big(x^2,\tilde{x}^3\big)$ such that
\begin{gather*}
\phi\big(x^2,\tilde{x}^3\big)=\Phi_{\tilde{x}^3},
\qquad
\psi\big(x^2,\tilde{x}^3\big)=\Phi_{x^2}.
\end{gather*}
Then equation~\eqref{case-2-3-eta1-in-coords} implies that
\begin{gather*}
x^1=\tilde{x}^1+\Phi\big(x^2,\tilde{x}^3\big)-\tilde{x}^3\Phi_{\tilde{x}^3}\big(x^2,\tilde{x}^3\big).
\end{gather*}
The local coordinate transformations which preserve the expression for $\eta^1$ in~\eqref{case-2-3-coframing} are
def\/ined implicitly by{\samepage
\begin{gather}
x^1=\tilde{x}^1+\Phi\big(x^2,\tilde{x}^3\big)-\tilde{x}^3\Phi_{\tilde{x}^3}\big(x^2,\tilde{x}^3\big),
\qquad
\tilde{x}^2=\Phi_{\tilde{x}^3}\big(x^2,\tilde{x}^3\big),
\qquad
x^3=\Phi_{x^2}\big(x^2,\tilde{x}^3\big),
\label{case-2-3-local-coord-trans}
\end{gather}
where $\Phi\big(x^2,\tilde{x}^3\big)$ is an arbitrary smooth function of two variables with $\Phi_{x^2\tilde{x}^3}\neq0$.}

Next we will compute how the function $H\big(x^1,x^2,x^3\big)$ transforms under a~coordinate transformation of the
form~\eqref{case-2-3-local-coord-trans}.
(When we consider the implications of homogeneity, it will turn out that $J$ can be expressed in terms of $H$ and its
derivatives; thus there is no need to explicitly compute the ef\/fects of the
transformation~\eqref{case-2-3-local-coord-trans} on $J$.) Consider the expression for $\eta^3$
in~\eqref{case-2-3-coframing}.
We must have
\begin{gather}
\eta^3=\frac{1}{\sqrt{\epsilon H_{x^1}(x)}}\big(H(x)\,dx^2-dx^3\big)=
\frac{1}{\sqrt{\epsilon\tilde{H}_{\tilde{x}^1}(\tilde{x})}}\big(\tilde{H}(\tilde{x})\,d\tilde{x}^2-d\tilde{x}^3\big).
\label{case-2-3-eta3-in-both-coords}
\end{gather}
From~\eqref{case-2-3-local-coord-trans}, we have
\begin{gather*}
d\tilde{x}^2=\Phi_{x^2\tilde{x}^3}\,dx^2+\Phi_{\tilde{x}^3\tilde{x}^3}\,d\tilde{x}^3,
\qquad
dx^3=\Phi_{x^2x^2}\,dx^2+\Phi_{x^2\tilde{x}^3}\,d\tilde{x}^3.
\end{gather*}
Substituting these expressions into~\eqref{case-2-3-eta3-in-both-coords} yields
\begin{gather}
\frac{1}{\sqrt{\epsilon H_{x^1}(x)}}\big((H(x)-\Phi_{x^2x^2})\,dx^2-\Phi_{x^2\tilde{x}^3}\,d\tilde{x}^3\big)
\nonumber
\\
\qquad
=\frac{1}{\sqrt{\epsilon\tilde{H}_{\tilde{x}^1}(\tilde{x})}}
\big(\tilde{H}(\tilde{x})\Phi_{x^2\tilde{x}^3}\,dx^2
+\big(\tilde{H}(\tilde{x})\Phi_{\tilde{x}^3\tilde{x}^3}-1\big) d\tilde{x}^3\big).
\label{case-2-3-eta3-computations1}
\end{gather}
Equating the ratios of the coef\/f\/icients of $dx^2$ and $d\tilde{x}^3$ on both sides
of~\eqref{case-2-3-eta3-computations1} yields
\begin{gather*}
\frac{(H(x)-\Phi_{x^2x^2})}{-\Phi_{x^2\tilde{x}^3}}=
\frac{\tilde{H}(\tilde{x})\Phi_{x^2\tilde{x}^3}}{(\tilde{H}(\tilde{x})\Phi_{\tilde{x}^3\tilde{x}^3}-1)},
\end{gather*}
which implies that
\begin{gather}
H\big(x^1,x^2,x^3\big)=
\frac{\big((\Phi_{x^2\tilde{x}^3})^2-\Phi_{x^2x^2}\Phi_{\tilde{x}^3\tilde{x}^3}\big)
\tilde{H}\big(\tilde{x}^1,\tilde{x}^2,\tilde{x}^3\big)
+\Phi_{x^2x^2}}{1-\Phi_{\tilde{x}^3\tilde{x}^3}\tilde{H}(\tilde{x}^1,\tilde{x}^2,\tilde{x}^3)}.
\label{case-2-3-new-H}
\end{gather}

Now suppose that the structure is homogeneous.
Unlike the previous cases, the assumption of homogeneity will imply some relations among the constants appearing in the
structure equations~\eqref{case-2-3-structure-eqs}.
In the homogeneous case, the functions $T^i_{jk}$ are all constant, and dif\/ferentiating
equations~\eqref{case-2-3-structure-eqs} implies that
\begin{gather*}
0=d\big(d\eta^1\big)=\big(T^1_{23}T^3_{13}-T^1_{13}T^3_{23}\big)\,\eta^1\wedge\eta^2\wedge\eta^3,
\\
0=d\big(d\eta^2\big)=\big(T^2_{23}T^3_{13}-T^2_{13}T^3_{23}\big)\,\eta^1\wedge\eta^2\wedge\eta^3,
\\
0=d\big(d\eta^3\big)=-\big(T^1_{13}+T^2_{23}\big)\,\eta^1\wedge\eta^2\wedge\eta^3.
\end{gather*}
The f\/irst two equations imply that the vectors
\begin{gather*}
\big[T^1_{13}
\;
T^1_{23}\big],
\qquad
\big[T^2_{13}
\;
T^2_{23}\big],
\qquad
\big[T^3_{13}
\;
T^3_{23}\big]
\end{gather*}
are all scalar multiples of each other unless $T^3_{13}=T^3_{23}=0$, while the third equation implies that
\begin{gather*}
T^2_{23}=-T^1_{13}.
\end{gather*}
In most of the computations that follow, these relations will be self-evident; however, at one point they will have
implications for the function~$H$.

The structure equation for $d\eta^1$ is
\begin{gather*}
d\eta^1=-J\sqrt{\epsilon H_{x^1}} \eta^1\wedge\eta^3-\epsilon \eta^2\wedge\eta^3.
\end{gather*}
Therefore, we must have
\begin{gather*}
J=\frac{c_1}{\sqrt{\epsilon H_{x^1}}}
\end{gather*}
for some constant $c_1$, so that the equation for $d\eta^1$ becomes
\begin{gather*}
d\eta^1=-c_1 \eta^1\wedge\eta^3-\epsilon \eta^2\wedge\eta^3.
\end{gather*}

Now the equation for $d\eta^3$ reduces to
\begin{gather*}
d\eta^3=\eta^1\wedge\eta^2\\
\phantom{d\eta^3=}
-\frac{1}{2H_{x^1}\sqrt{\epsilon H_{x^1}}}
\left(H_{x^1x^2}+x^3H_{x^1x^1}+H H_{x^1x^3}-2H_{x^1}H_{x^3}\right)
\left(c_1\eta^1\wedge\eta^3+\epsilon \eta^2\wedge\eta^3\right).
\end{gather*}
Therefore,
\begin{gather}
\frac{H_{x^1x^2}+x^3H_{x^1x^1}+H H_{x^1x^3}-2H_{x^1}H_{x^3}}{H_{x^1}\sqrt{\epsilon H_{x^1}}}=-2c_2
\label{case-2-3-pde1-for-H}
\end{gather}
for some constant $c_2$.
Substituting the derivative of~\eqref{case-2-3-pde1-for-H} with respect to $x^1$ into the equation for $d\eta^2$ yields
\begin{gather*}
d\eta^2=
\left(\frac{3}{4}\left(\frac{H_{x^1x^1}}{H_{x^1}}\right)^2-\frac{1}{2}\frac{H_{x^1x^1x^1}}{H_{x^1}}+\frac{c_1^2}{\epsilon}\right)\,\eta^1\wedge\eta^3+c_1\,\eta^2\wedge\eta^3.
\end{gather*}
Observe that:
\begin{itemize}\itemsep=0pt
\item The coef\/f\/icient of $\eta^2\wedge\eta^3$ in $d\eta^2$ is equal to minus the coef\/f\/icient of
$\eta^1\wedge\eta^3$ in $d\eta^1$, as we previously observed that it must be.
\item If $c_2\neq0$, then the ratio of the $\eta^1\wedge\eta^3$ and $\eta^2\wedge\eta^3$ coef\/f\/icients in $d\eta^2$
must be equal to $\frac{c_1}{\epsilon}$ (which is the ratio of these coef\/f\/icients in $d\eta^1$), and hence the
$\eta^1\wedge\eta^3$ coef\/f\/icient in $d\eta^2$ must be equal to $\frac{c_1^2}{\epsilon}$.
\end{itemize}
Therefore, if $c_2\neq0$, then $H$ satisf\/ies the PDE
\begin{gather}
\frac{3}{4}\left(\frac{H_{x^1x^1}}{H_{x^1}}\right)^2-\frac{1}{2}\frac{H_{x^1x^1x^1}}{H_{x^1}}=0.
\label{case-2-3-pde2-for-H}
\end{gather}
The solutions of~\eqref{case-2-3-pde2-for-H} are precisely the linear fractional transformations in the $x^1$ variable,
and so we must have
\begin{gather}
H\big(x^1,x^2,x^3\big)=\frac{F_1\big(x^2,x^3\big)x^1+F_0\big(x^2,x^3\big)}{G_1\big(x^2,x^3\big)x^1+G_0\big(x^2,x^3\big)}
\label{case-2-3-reduced-H-1}
\end{gather}
for some functions $F_0\big(x^2,x^3\big)$, $F_1\big(x^2,x^3\big)$, $G_0\big(x^2,x^3\big)$, $G_1\big(x^2,x^3\big)$.

By contrast, if $c_2=0$, then the vectors $[T^1_{13}\;T^1_{23}]$, $[T^2_{13}\;T^2_{23}]$
are no longer required to be linearly independent, and so the Schwarzian derivative of $H$ with respect to
$x^1$ appearing in equation~\eqref{case-2-3-pde2-for-H} is only required to be constant, but not necessarily equal to
zero.
There are two possibilities, depending on the sign:
\begin{itemize}
\itemsep=0pt \item If
$\frac{3}{4}\left(\frac{H_{x^1x^1}}{H_{x^1}}\right)^2-\frac{1}{2}\frac{H_{x^1x^1x^1}}{H_{x^1}}=-c_3^2$ for $c_3>0$, then
\begin{gather}
H\big(x^1,x^2,x^3\big)=F_1\big(x^2,x^3\big)\tan\left(c_3x^1+F_0\big(x^2,x^3\big)\right)+F_2\big(x^2,x^3\big)
\label{case-2-3-H-tan-case}
\end{gather}
for some functions $F_0\big(x^2,x^3\big)$, $F_1\big(x^2,x^3\big)$, $F_2\big(x^2,x^3\big)$ with $F_1\big(x^2,x^3\big)\neq0$.
\item If $\frac{3}{4}\left(\frac{H_{x^1x^1}}{H_{x^1}}\right)^2-\frac{1}{2}\frac{H_{x^1x^1x^1}}{H_{x^1}}=c_3^2$ for
$c_3>0$, then
\begin{gather}
H\big(x^1,x^2,x^3\big)=F_1\big(x^2,x^3\big)\tanh\left(c_3x^1+F_0\big(x^2,x^3\big)\right)+F_2\big(x^2,x^3\big)
\label{case-2-3-H-tanh-case}
\end{gather}
for some functions $F_0\big(x^2,x^3\big)$, $F_1\big(x^2,x^3\big)$, $F_2\big(x^2,x^3\big)$ with $F_1\big(x^2,x^3\big)\neq0$.
\end{itemize}
We consider each of these cases separately.

\subsection[c_2\neq0]{$\boldsymbol{c_2\neq0}$}

In this case, $H\big(x^1,x^2,x^3\big)$ is given by~\eqref{case-2-3-reduced-H-1}.
Now we compute how the function~\eqref{case-2-3-reduced-H-1} transforms under a~local coordinate transformation of the
form~\eqref{case-2-3-local-coord-trans}.
\begin{Lemma}
There exists a~local coordinate transformation of the form~\eqref{case-2-3-local-coord-trans} such that the function
$\tilde{H}(\tilde{x}^1,\tilde{x}^2,\tilde{x}^3)$ is linear in $\tilde{x}^1$, i.e.,
\begin{gather}
\tilde{H}\big(\tilde{x}^1,\tilde{x}^2,\tilde{x}^3\big)=
\tilde{F}_1\big(\tilde{x}^2,\tilde{x}^3\big)\tilde{x}^1+\tilde{F}_0\big(\tilde{x}^2,\tilde{x}^3\big),
\label{case-2-3-H-linear-form}
\end{gather}
with $\tilde{F}_1\neq0$.
\end{Lemma}
\begin{proof}
Equation~\eqref{case-2-3-new-H} can be written as
\begin{gather}
\tilde{H}(\tilde{x}^1,\tilde{x}^2,\tilde{x}^3)=
\frac{H\big(x^1,x^2,x^3\big)-\Phi_{x^2x^2}}{\Phi_{\tilde{x}^3\tilde{x}^3}H\big(x^1,x^2,x^3\big)+\big((\Phi_{x^2\tilde{x}^3})^2-
\Phi_{x^2x^2}\Phi_{\tilde{x}^3\tilde{x}^3}\big)}.
\label{case-2-3-old-H}
\end{gather}
Substituting~\eqref{case-2-3-reduced-H-1} into this equation yields
\begin{gather*}
\begin{split}
& \frac{\tilde{F}_1\tilde{x}^1+\tilde{F}_0}{\tilde{G}_1\tilde{x}^1+\tilde{G}_0}=
\frac{\left(\frac{F_1x^1+F_0}{G_1x^1+G_0}\right)-\Phi_{x^2x^2}}{\Phi_{\tilde{x}^3\tilde{x}^3}\left(\frac{F_1x^1+F_0}{G_1x^1+G_0}\right)
+\left((\Phi_{x^2\tilde{x}^3})^2-\Phi_{x^2x^2}\Phi_{\tilde{x}^3\tilde{x}^3}\right)}
\\
& \qquad
=\frac{\left(F_1x^1+F_0\right)-\Phi_{x^2x^2}\left(G_1x^1+G_0\right)}{\Phi_{\tilde{x}^3\tilde{x}^3}\left(F_1x^1+F_0\right)
+\left((\Phi_{x^2\tilde{x}^3})^2-\Phi_{x^2x^2}\Phi_{\tilde{x}^3\tilde{x}^3}\right)\left(G_1x^1+G_0\right)}
\\
& \qquad
=\frac{\left[F_1-\Phi_{x^2x^2}G_1\right]x^1+\left[F_0-\Phi_{x^2x^2}G_0\right]}
{\left[\Phi_{\tilde{x}^3\tilde{x}^3}F_1
\!+\!\left((\Phi_{x^2\tilde{x}^3})^2\!-\!\Phi_{x^2x^2}\Phi_{\tilde{x}^3\tilde{x}^3}\right)G_1\right]x^1
+\left[\Phi_{\tilde{x}^3\tilde{x}^3}F_0
\!+\!\left((\Phi_{x^2\tilde{x}^3})^2\!-\!\Phi_{x^2x^2}\Phi_{\tilde{x}^3\tilde{x}^3}\right)G_0\right]}.
\end{split}
\end{gather*}
The coef\/f\/icients of $\tilde{x}^1$ on the left-hand side of this equation are the same as the coef\/f\/icients of~$x^1$ on the right-hand side, so the condition that $\tilde{G}_1=0$ is equivalent to
\begin{gather*}
0=\Phi_{\tilde{x}^3\tilde{x}^3}\big(x^2,\tilde{x}^3\big)F_1\big(x^2,x^3\big)
+\left(\big(\Phi_{x^2\tilde{x}^3}\big(x^2,\tilde{x}^3\big)\big)^2\!-\Phi_{x^2x^2}\big(x^2,\tilde{x}^3\big)\Phi_{\tilde{x}^3\tilde{x}^3}\big(x^2,\tilde{x}^3\big)\right)G_1\big(x^2,x^3\big)
\\
\phantom{0}
=\Phi_{\tilde{x}^3\tilde{x}^3}\big(x^2,\tilde{x}^3\big)F_1\big(x^2,\Phi_{x^2}\big(x^2,\tilde{x}^3\big)\big)
\\
\phantom{0=}
+\left(\big(\Phi_{x^2\tilde{x}^3}\big(x^2,\tilde{x}^3\big)\big)^2
-\Phi_{x^2x^2}\big(x^2,\tilde{x}^3\big)\Phi_{\tilde{x}^3\tilde{x}^3}\big(x^2,\tilde{x}^3\big)\right)
G_1\big(x^2,\Phi_{x^2}\big(x^2,\tilde{x}^3\big)\big).
\end{gather*}
Any solution $\Phi\big(x^2,\tilde{x}^3\big)$ of this equation will induce a~local coordinate transformation for which
$\tilde{H}(\tilde{x}^1,\tilde{x}^2,\tilde{x}^3)$ has the form~\eqref{case-2-3-H-linear-form}, as desired.
Note that $\tilde{F}_1=\tilde{H}_{x^1}\neq0$, and hence $\tilde{F}_1$ must have the same sign as $\epsilon$.
\end{proof}

Local coordinates for which $H$ has the form~\eqref{case-2-3-H-linear-form} are determined up to transformations of the
form~\eqref{case-2-3-local-coord-trans} with $\Phi_{\tilde{x}^3\tilde{x}^3}=0$, i.e.,
\begin{gather*}
\Phi\big(x^2,\tilde{x}^3\big)=\Phi_1\big(x^2\big)\tilde{x}^3+\Phi_0\big(x^2\big).
\end{gather*}
With $\Phi$ as above, the local coordinate transformation~\eqref{case-2-3-local-coord-trans} reduces to
\begin{gather}\label{case-2-3-local-coord-trans-2}
x^1=\tilde{x}^1+\Phi_0\big(x^2\big),
\qquad
\tilde{x}^2=\Phi_1\big(x^2\big),
\qquad
x^3=\Phi_0'\big(x^2\big)+\Phi_1'\big(x^2\big)\tilde{x}^3.
\end{gather}

With the assumption that $H$ has the form
\begin{gather*}
H\big(x^1,x^2,x^3\big)=F_1\big(x^2,x^3\big)x^1+F_0\big(x^2,x^3\big),
\end{gather*}
dif\/ferentiating equation~\eqref{case-2-3-pde1-for-H} with respect to $x^1$ yields
\begin{gather*}
\frac{(F_1)_{x^3}}{\sqrt{\epsilon F_1}}=0.
\end{gather*}
Therefore,
\begin{gather*}
F_1\big(x^2,x^3\big)=F_1\big(x^2\big).
\end{gather*}
Now equation~\eqref{case-2-3-new-H} reduces to
\begin{gather*}
F_1\big(x^2\big)x^1+F_0\big(x^2,x^3\big)
=\left(\Phi_1'\big(x^2\big)\right)^2
\big(\tilde{F}_1\big(\tilde{x}^2\big)\tilde{x}^1+\tilde{F}_0\big(\tilde{x}^2,\tilde{x}^3\big)\big)
+\Phi_0''\big(x^2\big)+\Phi_1''\big(x^2\big)\tilde{x}^3,
\end{gather*}
which, taking~\eqref{case-2-3-local-coord-trans-2} into account, becomes
\begin{gather}
F_1\big(x^2\big)\tilde{x}^1+\left(F_1\big(x^2\big)\Phi_0\big(x^2\big)
+F_0\big(x^2,\Phi_0'\big(x^2\big)+\Phi_1'\big(x^2\big)\tilde{x}^3\big)\!\right)
\nonumber
\\
\qquad
=\left(\Phi_1'\big(x^2\big)\!\right)^2\tilde{F}_1\big(\Phi_1\big(x^2\big)\!\big)\tilde{x}^1+
\left(\Phi_1'\big(x^2\big)\!\right)^2\tilde{F}_0\big(\Phi_1\big(x^2\big),\tilde{x}^3\big)
+\Phi_0''\big(x^2\big)+\Phi_1''\big(x^2\big)\tilde{x}^3.
\label{case-2-3-new-H-2}
\end{gather}
Equating the coef\/f\/icients of $\tilde{x}^1$ on both sides yields
\begin{gather*}
F_1\big(x^2\big)=\left(\Phi_1'\big(x^2\big)\right)^2\tilde{F}_1\big(\Phi_1\big(x^2\big)\big).
\end{gather*}
Thus any solution $\Phi_1\big(x^2\big)$ of the equation
\begin{gather*}
\Phi_1'\big(x^2\big)=\sqrt{\epsilon F_1\big(x^2\big)}
\end{gather*}
will induce a~local coordinate transformation for which
\begin{gather*}
\tilde{F}_1\big(\tilde{x}^2\big)=\epsilon=\pm1.
\end{gather*}
Local coordinates for which $F_1\big(x^2\big)=\epsilon$ are determined up to transformations of the
form~\eqref{case-2-3-local-coord-trans-2} with
\begin{gather*}
\Phi_1'\big(x^2\big)=\pm1;
\end{gather*}
for simplicity, we will assume that $\Phi_1'\big(x^2\big)=1$.
Then
\begin{gather*}
\Phi\big(x^2,\tilde{x}^3\big)=x^2\tilde{x}^3+a\tilde{x}^3+\Phi_0\big(x^2\big)
\end{gather*}
for some constant $a$.
With $\Phi$ as above, the local coordinate transformation~\eqref{case-2-3-local-coord-trans} reduces to
\begin{gather}
\label{case-2-3-local-coord-trans-3}
x^1=\tilde{x}^1+\Phi_0\big(x^2\big),
\qquad
\tilde{x}^2=x^2+a,
\qquad
x^3=\tilde{x}^3+\Phi_0'\big(x^2\big).
\end{gather}

Now equation~\eqref{case-2-3-pde1-for-H} takes the form
\begin{gather*}
(F_0)_{x^3}=\epsilon c_3.
\end{gather*}
Therefore,
\begin{gather*}
F_0\big(x^2,x^3\big)=\epsilon c_3x^3+F_{2}\big(x^2\big).
\end{gather*}
Now equation~\eqref{case-2-3-new-H-2} reduces to
\begin{gather*}
\epsilon\Phi_0\big(x^2\big)+\epsilon c_3\Phi_0'\big(x^2\big)+F_2\big(x^2\big)=\tilde{F}_2(x^2+a)+\Phi_0''\big(x^2\big).
\end{gather*}
Thus any solution $\Phi_0\big(x^2\big)$ of the equation
\begin{gather*}
\epsilon\Phi_0\big(x^2\big)+\epsilon c_3\Phi_0'\big(x^2\big)-\Phi_0''\big(x^2\big)=-F_2\big(x^2\big)
\end{gather*}
will induce a~local coordinate transformation for which
\begin{gather*}
\tilde{F}_2\big(\tilde{x}^2\big)=0.
\end{gather*}
Local coordinates for which $F_2\big(x^2\big)=0$ are determined up to transformations of the
form~\eqref{case-2-3-local-coord-trans-3} with
\begin{gather*}
\epsilon\Phi_0\big(x^2\big)+\epsilon c_3\Phi_0'\big(x^2\big)-\Phi_0''\big(x^2\big)=0,
\end{gather*}
i.e.,
\begin{gather*}
\Phi_0\big(x^2\big)=b_1e^{r_1x^2}+b_2e^{r_2x^2},
\end{gather*}
where $b_1$, $b_2$ are constants and
\begin{gather}
r_1=\frac{\epsilon c_3+\sqrt{c_3^2+4\epsilon}}{2},
\qquad
r_2=\frac{\epsilon c_3-\sqrt{c_3^2+4\epsilon}}{2}.
\label{case-2-3-define-r's}
\end{gather}
Note that if $\epsilon=1$, then $r_1$, $r_2$ are real and distinct; if $\epsilon=-1$, then $r_1$, $r_2$ may be real and
distinct, real and equal, or a~complex conjugate pair.

To summarize, we have constructed local coordinates for which
\begin{gather*}
J\big(x^1,x^2,x^3\big)=c_1,
\qquad
H\big(x^1,x^2,x^3\big)=\epsilon\big(x^1+c_3x^3\big).
\end{gather*}
These coordinates are determined up to transformations of the form
\begin{gather*}
x^1=\tilde{x}^1+b_1e^{r_1x^2}+b_2e^{r_2x^2},
\qquad
\tilde{x}^2=x^2+a,
\qquad
x^3=\tilde{x}^3+b_1r_1e^{r_1x^2}+b_2r_2e^{r_2x^2}.
\end{gather*}

\subsection[$c_2=0$]{$\boldsymbol{c_2=0}$}
We will only give the details of the analysis for the case~\eqref{case-2-3-H-tan-case}; the
case~\eqref{case-2-3-H-tanh-case} is similar.
First we compute how the function~\eqref{case-2-3-H-tan-case} transforms under a~local coordinate transformation of the
form~\eqref{case-2-3-local-coord-trans}.
\begin{Lemma}
There exists a~local coordinate transformation of the form~\eqref{case-2-3-local-coord-trans} such that \newline
$\tilde{F}_0\big(\tilde{x}^2,\tilde{x}^3\big)=0$, i.e.,
\begin{gather}
\tilde{H}\big(\tilde{x}^1,\tilde{x}^2,\tilde{x}^3\big)=
\tilde{F}_1\big(\tilde{x}^2,\tilde{x}^3\big)\tan(c_3\tilde{x}^1)+\tilde{F}_2\big(\tilde{x}^2,\tilde{x}^3\big),
\label{case-2-3-H-tan-linear-form}
\end{gather}
with $\tilde{F}_1\neq0$.

\end{Lemma}
\begin{proof}
Substituting~\eqref{case-2-3-H-tan-case} into the expression~\eqref{case-2-3-old-H} for
$\tilde{H}(\tilde{x}^1,\tilde{x}^2,\tilde{x}^3)$ yields
\begin{gather}
\tilde{F}_1\tan\big(c_3\tilde{x}^1+\tilde{F}_0\big)+\tilde{F}_2
\nonumber
\\
\qquad
=\frac{\left(F_1\tan(c_3x^1+F_0)+F_2\right)-\Phi_{x^2x^2}}
{\Phi_{\tilde{x}^3\tilde{x}^3}\left(F_1\tan(c_3x^1+F_0)+F_2\right)
+\left((\Phi_{x^2\tilde{x}^3})^2-\Phi_{x^2x^2}\Phi_{\tilde{x}^3\tilde{x}^3}\right)}
\nonumber
\\
\qquad
=\frac{F_1\sin(c_3x^1+F_0)+(F_2-\Phi_{x^2x^2})\cos(c_3x^1+F_0)}
{\Phi_{\tilde{x}^3\tilde{x}^3}F_1\sin(c_3x^1\!+\!F_0)
\!+\!\left(\Phi_{\tilde{x}^3\tilde{x}^3}F_2\!+\!
(\Phi_{x^2\tilde{x}^3})^2\!-\!\Phi_{x^2x^2}\Phi_{\tilde{x}^3\tilde{x}^3}\right)
\cos(c_3x^1\!+\!F_0)}.
\label{case-2-3-2-new-H-yuck1}
\end{gather}
Now, def\/ine functions $R\big(x^2,\tilde{x}^3\big)$, $\Theta\big(x^2,\tilde{x}^3\big)$ by the conditions that
\begin{gather*}
\Phi_{\tilde{x}^3\tilde{x}^3}F_1=R\sin\Theta,
\qquad
\Phi_{\tilde{x}^3\tilde{x}^3}F_2+(\Phi_{x^2\tilde{x}^3})^2-\Phi_{x^2x^2}\Phi_{\tilde{x}^3\tilde{x}^3}=R\cos\Theta;
\end{gather*}
in particular, we have
\begin{gather*}
\Theta=
\tan^{-1}\left(\frac{\Phi_{\tilde{x}^3\tilde{x}^3}F_1}{\Phi_{\tilde{x}^3\tilde{x}^3}F_2+(\Phi_{x^2\tilde{x}^3})^2-\Phi_{x^2x^2}\Phi_{\tilde{x}^3\tilde{x}^3}}\right).
\end{gather*}
Then the denominator of the right-hand side of~\eqref{case-2-3-2-new-H-yuck1} can be written as
\begin{gather*}
R\cos\big(c_3x^1+F_0-\Theta\big)
=R\cos\big(c_3\big(\tilde{x}^1+\Phi-\tilde{x}^3\Phi_{\tilde{x}^3}\big)+F_0-\Theta\big).
\end{gather*}
Therefore, $\tilde{H}(\tilde{x}^1,\tilde{x}^2,\tilde{x}^3)$ is a~linear function of the quantity
\begin{gather*}
\tan\big(c_3\big(\tilde{x}^1+\Phi-\tilde{x}^3\Phi_{\tilde{x}^3}\big)+F_0-\Theta\big),
\end{gather*}
which implies that
\begin{gather*}
\tilde{F}_0=c_3\big(\Phi-\tilde{x}^3\Phi_{\tilde{x}^3}\big)+F_0-\Theta.
\end{gather*}
Keeping in mind that $\Theta$ is a~second-order dif\/ferential operator in $\Phi$, the condition $\tilde{F}_0=0$ is
a~second-order PDE for the function $\Phi\big(x^2,\tilde{x}^3\big)$.
Any solution $\Phi\big(x^2,\tilde{x}^3\big)$ of this equation will induce a~local coordinate transformation for which
$\tilde{H}(\tilde{x}^1,\tilde{x}^2,\tilde{x}^3)$ has the form~\eqref{case-2-3-H-tan-linear-form}, as desired.
\end{proof}

Local coordinates for which $H$ has the form~\eqref{case-2-3-H-tan-linear-form} are determined up to transformations of
the form~\eqref{case-2-3-local-coord-trans} with $\Phi\big(x^2,\tilde{x}^3\big)$ a~solution of the PDE
\begin{gather}
c_3\big(\Phi-\tilde{x}^3\Phi_{\tilde{x}^3}\big)-\Theta=0.
\label{case-2-3-2-fugly-pde}
\end{gather}
Unfortunately we cannot explicitly write down the general solution to this PDE; however, a~subset of the solutions is
given by the family{\samepage
\begin{gather*}
\Phi\big(x^2,\tilde{x}^3\big)=\tilde{x}^3\Phi_0\big(x^2\big),
\end{gather*}
where $\Phi_0$ is an arbitrary function of $x^2$ with $\Phi_0'\big(x^2\big)\neq0$.}

With the assumption that $H$ has the form~\eqref{case-2-3-H-tan-linear-form}, equation~\eqref{case-2-3-pde1-for-H}
becomes (recalling that $c_2=0$)
\begin{gather*}
F_1\big((F_1)_{{x}^3}-2c_3x^3\big)\tan\big(c_3x^1\big)+2F_1(F_2)_{{x}^3}-F_2(F_1)_{{x}^3}-(F_1)_{{x}^2}=0.
\end{gather*}
Therefore, since $F_1\neq0$,
\begin{gather}
\label{case-3-2-F1-F2-pdes}
(F_1)_{{x}^3}-2c_3x^3=0,
\qquad
2F_1(F_2)_{{x}^3}-F_2(F_1)_{{x}^3}-(F_1)_{{x}^2}=0.
\end{gather}
The f\/irst equation implies that
\begin{gather*}
F_1\big(x^2,x^3\big)=c_3\big(x^3\big)^2+F_{10}\big(x^2\big)
\end{gather*}
for some function $F_{10}\big(x^2\big)$.

Computations similar to those in the previous case show that, under a~local coordinate
transformation~\eqref{case-2-3-local-coord-trans} with $\Phi=\tilde{x}^3\Phi_0\big(x^2\big)$, we have
\begin{gather*}
\tilde{F}_{10}\big(\tilde{x}^2\big)=\frac{F_{10}\big(x^2\big)}{\Phi_0'\big(x^2\big)^2}.
\end{gather*}
Thus any solution $\Phi_0\big(x^2\big)$ of the equation
\begin{gather*}
\Phi_0'\big(x^2\big)=\sqrt{\frac{1}{c_3}\epsilon F_{10}\big(x^2\big)}
\end{gather*}
will induce a~local coordinate transformation for which either $\tilde{F}_{10}\big(\tilde{x}^2\big)=0$ or
$\tilde{F}_{10}\big(\tilde{x}^2\big)=\pm c_3$.
Denote this constant by $c_4$, so that we now have
\begin{gather*}
\tilde{F}_1\big(\tilde{x}^2,\tilde{x}^3\big)=c_3\big(\tilde{x}^3\big)^2+c_4.
\end{gather*}

Finally, the second equation in~\eqref{case-3-2-F1-F2-pdes} becomes
\begin{gather*}
2\big(c_3\big(x^3\big)^2+c_4\big)(F_2)_{{x}^3}-2c_3x^3F_2=0,
\end{gather*}
which implies that
\begin{gather*}
F_2\big(x^2,x^3\big)=F_{20}\big(x^2\big)\sqrt{c_3(x^3)^2+c_4}
\end{gather*}
for some function $F_{20}\big(x^2\big)$.
We conjecture that the remaining solutions of~\eqref{case-2-3-2-fugly-pde} can be used to normalize the function
$F_{20}\big(x^2\big)$, but unfortunately we have been unable to complete this step in the analysis.

To summarize, we have constructed local coordinates for which
\begin{gather*}
J\big(x^1,x^2,x^3\big)=\frac{c_1\cos(c_3x^1)}{\sqrt{c_3(c_2(x^3)^2+c_4)}},
\\
H\big(x^1,x^2,x^3\big)=\Big(c_3\big(x^3\big)^2+c_4\Big)\tan\big(c_3x^1\big)+F_{20}\big(x^2\big)\sqrt{c_3(x^3)^2+c_4}.
\end{gather*}

\subsection*{Acknowledgements} This research was supported in part by NSF grants DMS-0908456 and DMS-1206272.
We would like to thank the referees for many helpful suggestions, which signif\/icantly improved the organization and
exposition of this paper.

\pdfbookmark[1]{References}{ref}
\LastPageEnding


\begin{thebibliography}{99} \footnotesize\itemsep=0pt

\bibitem{CMW09}
Clelland J.N., Moseley C.G., Wilkens G.R., Geometry of control-affine systems,
  \href{http://dx.doi.org/10.3842/SIGMA.2009.095}{\textit{SIGMA}} \textbf{5} (2009), 095, 28~pages, \href{http://arxiv.org/abs/0903.4932}{arXiv:0903.4932}.

\bibitem{Gardner89}
Gardner R.B., The method of equivalence and its applications, \href{http://dx.doi.org/10.1137/1.9781611970135}{\textit{CBMS-NSF
  Regional Conference Series in Applied Mathematics}}, Vol.~58, Society for
  Industrial and Applied Mathematics (SIAM), Philadelphia, PA, 1989.

\bibitem{Jurdjevic97}
Jurdjevic V., Geometric control theory, \textit{Cambridge Studies in Advanced
  Mathematics}, Vol.~52, Cambridge University Press, Cambridge, 1997.

\bibitem{LS95}
Liu W., Sussman H.J., Shortest paths for sub-{R}iemannian metrics on rank-two
  distributions, \textit{Mem. Amer. Math. Soc.} \textbf{118} (1995), no.~564,
  x+104~pages.

\bibitem{Montgomery02}
Montgomery R., A tour of subriemannian geometries, their geodesics and
  applications, \textit{Mathematical Surveys and Monographs}, Vol.~91, Amer.
  Math. Soc., Providence, RI, 2002.

\bibitem{Moseley01}
Moseley C.G., The geometry of sub-{R}iemannian {E}ngel manifolds, Ph.D.
  thesis, University of North Caro\-li\-na at Chapel Hill, 2001.

\end{thebibliography}
\end{document}